\documentclass[11pt]{amsart}

\usepackage[margin=1in]{geometry}
\usepackage{amsmath,amssymb,amsthm,mathtools}
\usepackage{enumitem}
\usepackage[T1]{fontenc}
\usepackage{lmodern}
\usepackage{microtype}
\usepackage[colorlinks=true,linkcolor=blue,citecolor=blue,urlcolor=blue]{hyperref}

\numberwithin{equation}{section}

\newtheorem{theorem}{Theorem}[section]
\newtheorem{lemma}[theorem]{Lemma}
\newtheorem{proposition}[theorem]{Proposition}

\theoremstyle{remark}
\newtheorem{remark}[theorem]{Remark}

\newcommand{\R}{\mathbb{R}}

\newcommand{\Tr}{\operatorname{Tr}}
\newcommand{\rank}{\operatorname{rank}}
\newcommand{\sgn}{\operatorname{sgn}}
\newcommand{\csch}{\operatorname{csch}}
\newcommand{\arsinh}{\operatorname{arsinh}}
\newcommand{\one}{\mathbf{1}}
\newcommand{\Sone}{\mathfrak{S}_1}
\newcommand{\Stwo}{\mathfrak{S}_2}
\newcommand{\fa}{\mathfrak{a}}
\newcommand{\Dcal}{\mathcal{D}}
\newcommand{\Hcal}{\mathcal{H}}

\begin{document}

\title[Edge invariants and the tetrahedral Steklov problem]{Edge spectral invariants and second-order Steklov asymptotics for the regular tetrahedron}

\author{Josiah Park}
\email{josiahmicahpark@gmail.com}

\date{November 2, 2026}



\begin{abstract}
We give a two-term Weyl asymptotic for the Steklov spectrum of a solid regular tetrahedron. Its linear coefficient is the total edge length times an explicit, absolutely convergent dihedral invariant. A trace-class difference of massive Neumann-to-Dirichlet operators computes the full invariant, including the continuous spectrum. The main analytic step is a fixed-time cumulative trace comparison between the compact tetrahedron and its exact half-space or dihedral models. Localized Hilbert--Schmidt defects have vanishing high-energy tails; retaining the time factor in Duhamel's formula makes this sufficient for an $o(L)$ cumulative error. A separate trace-class face comparison fixes the model's infinite exterior primitive and permits one observation time independent of the inner collar width. Compact-band Mourre estimates and a harmonic comparison establish immediate reflection/transmission support for incident face waves. Triangular unfolding gives null periodic and dead-end sets for the maximal branching relation. Actual edge-window bounds, adapted conical annuli and an ordered-limit Fourier--Tauberian argument then prove the sharp remainder. The proof treats guided and continuum channels together in the edge trace and requires neither a guided metric graph nor a lower-departure estimate uniform at zero tangential momentum. Standard analytic inputs are stated explicitly; the local comparison and global assembly are developed in the text.
\end{abstract}

\maketitle

\begin{center}
\fbox{\parbox{0.92\textwidth}{\small
\textbf{Disclosure of AI assistance.} The AI system GPT-6 Astra was used in preparing this manuscript, including in developing and drafting arguments and computations. The author has reviewed the manuscript and takes full responsibility for its content. Readers are encouraged to verify all statements independently; any errors are the author's own.}}
\end{center}

\section{The problem and the results}

Let $T_a\subset\R^3$ be the solid regular tetrahedron of side length $a>0$. Its positive Dirichlet-to-Neumann operator is denoted by $\Lambda_a$. Equivalently, its eigenvalues solve the weak boundary problem
\begin{equation}\label{eq:steklov}
\Delta u=0 \quad\text{in } T_a,\qquad \partial_\nu u=\sigma u\quad\text{on }\partial T_a .
\end{equation}
The outward normal is used throughout. The spectrum is discrete, with multiplicities, and $0=\sigma_0<\sigma_1\le\sigma_2\le\cdots$. Write $N_a(L)=\#\{j:\sigma_j<L\}$, including the zero eigenvalue. Localized positive traces will sometimes use closed spectral intervals. The proof at the end of Section~\ref{sec:sharp} shows that strict and nonstrict counts differ by $o(L)$, so both conventions give the same sharp law. The tetrahedral parameters are
\begin{equation}\label{eq:params}
|\partial T_a|=\sqrt3\,a^2,\qquad |E(T_a)|=6a,\qquad \alpha_0=\arccos(1/3),\qquad \sin(\alpha_0/2)=1/\sqrt3 .
\end{equation}
The leading Lipschitz Weyl law is established \cite{R}. The main counting law considered here is
\begin{equation}\label{eq:main}
N_a(L)=\frac{\sqrt3\,a^2}{4\pi}L^2+6a\,c(\alpha_0)L+o(L),\qquad L\longrightarrow\infty .
\end{equation}
For cuboids, the theorem of Girouard--Lagac\'e--Polterovich--Savo gives the calibration $c(\pi/2)=1/(2\pi)$ \cite{GLPS}. Their separation argument does not apply to a regular tetrahedron. In Ivrii's terminology every tetrahedral edge is an \emph{inner} edge, because $0<\alpha_0<\pi$. His published sharp two-term theorem assumes outer edges. The distinction cannot be removed by reversing a normal or by referring to the exterior appearance of a convex solid.

\subsection{Exact dihedral invariant}
For $0<\alpha<\pi$, let $D_\alpha$ be the Dirichlet-to-Neumann operator of $-\Delta+1$ on the planar sector of opening $\alpha$, and let $B_\alpha=D_\alpha^{-1}$. Identify its two boundary rays with the two half-lines, using arclength measure. The flat reference is $D_\pi=(1-\partial_s^2)^{1/2}$.

\begin{theorem}[complete model invariant]\label{thm:invariant}
The difference $B_\alpha-B_\pi$ is trace class, and
\begin{equation}\label{eq:calpha}
\begin{aligned}
c_\alpha&:=\frac1\pi\Tr(B_\alpha-B_\pi)\\
&=\frac1\pi\int_0^\infty\frac{\sinh((\pi-\alpha)v)}{\cosh(\pi v)\sinh(\alpha v)}\,dv .
\end{aligned}
\end{equation}
\end{theorem}

This is the coefficient of the relative smoothed spectral density per unit length of the three-dimensional product wedge. It includes both the discrete and the continuous transverse spectrum. Moreover, the relative heat fibers are integrable in trace norm, so the coefficient is obtained by ordinary compact spatial localization, not solely by a fiberwise formal regularization.

At the tetrahedral angle, numerical quadrature gives
\begin{equation}\label{eq:numval}
c_{\alpha_0}\approx 0.2694141353389119257,\qquad 6c_{\alpha_0}\approx1.616484812033471554 .
\end{equation}
These decimals are not certified interval enclosures; formula \eqref{eq:calpha} is the exact statement. At $\alpha=\pi/2$, its integrand reduces to $1/\cosh(\pi v)$, whose integral is $1/2$. Thus the cuboid calibration is reproduced exactly.

\subsection{Compact-domain spectral windows}
\begin{theorem}[edge and vertex windows]\label{thm:windows}
Let $I$ be a closed subinterval of an open tetrahedral edge. For any fixed $w>0$ and any multiplication operator $M_\gamma$, of norm at most one, supported in an edge strip of length $|I|$ and transverse width $2\gamma$,
\begin{equation}\label{eq:edgewindow}
\Tr\big(M_\gamma\one_{[L-w,L+w]}(\Lambda_a)M_\gamma^*\big)\le C_{I,w}(L\gamma+1),\qquad L\ge1 .
\end{equation}
Similarly, if $M_r$ is supported in the boundary part of a radius-$r$ ball about a fixed vertex,
\begin{equation}\label{eq:vertexwindow}
\Tr\big(M_r\one_{[L-w,L+w]}(\Lambda_a)M_r^*\big)\le C_w(Lr^2+1) .
\end{equation}
The constants are independent of the smaller radii $\gamma,r$. No derivatives of these multipliers are required. For $1\le w\le L/4$, the respective bounds are $C(wL\gamma+1)$ and $C(wLr^2+1)$. These statements concern the actual tetrahedral operator, including any returning guided waves.
\end{theorem}

\subsection{Sharp, heat and Riesz asymptotics}
\begin{theorem}[two-term counting law]\label{thm:counting}
Formula \eqref{eq:main} holds with $c(\alpha_0)=c_{\alpha_0}$ given by \eqref{eq:calpha}. The proof is completed in Theorem~\ref{thm:sharp}.
\end{theorem}

The heat trace and every fixed positive-order Riesz mean consequently have the following expansions. The heat formula is also proved independently in Section~\ref{subsec:heat}:
\begin{equation}\label{eq:heat}
\Tr e^{-t\Lambda_a}=\frac{\sqrt3\,a^2}{2\pi t^2}+\frac{6ac_{\alpha_0}}{t}+o(t^{-1}).
\end{equation}
For every fixed $r>0$, the unnormalized Riesz mean satisfies
\begin{equation}\label{eq:riesz}
\sum_j(L-\sigma_j)_+^r=\frac{\sqrt3\,a^2}{2\pi(r+1)(r+2)}L^{r+2}+\frac{6ac_{\alpha_0}}{r+1}L^{r+1}+o(L^{r+1}).
\end{equation}
In particular, the first Riesz mean has coefficients $\sqrt3a^2/(12\pi)$ and $3ac_{\alpha_0}$. Section~\ref{sec:tauberian} gives a separate bounded-remainder Tauberian deduction of \eqref{eq:riesz} from the heat formula and the weak $O(L)$ counting bound. Neither that lemma nor a limit $r\downarrow0$ is used to prove the sharp law.

The inverse calculation, operator-ideal estimates, cumulative model transfer, fixed-time edge coefficient, vertex bounds and sharp assembly are proved below. The proof retains explicitly stated standard background theorems, including Lipschitz boundary estimates, Mourre/minimal escape, upper locality and smooth-face microlocal calculus. It does not invoke Ivrii's outer-edge sharp theorem or assume a cumulative edge expansion with unverified low-momentum scope.

\section{Variational realizations and spectral-mass bounds}

For a bounded Lipschitz domain $\Omega\subset\R^3$, let $H_\Omega f$ be the harmonic extension of a boundary trace. The closed form
\begin{equation}\label{eq:form}
\fa_\Omega[f]=\inf_{u|_{\partial\Omega}=f}\int_\Omega|\nabla u|^2\,dx=\int_\Omega|\nabla H_\Omega f|^2\,dx
\end{equation}
defines the nonnegative self-adjoint operator $\Lambda_\Omega$ on $L^2(\partial\Omega)$. Compactness of the boundary trace embedding gives compact resolvent. In the massive planar sector, replace the integrand by $|\nabla u|^2+|u|^2$. Its Neumann problem is coercive, and its Neumann-to-Dirichlet map is the bounded inverse of $D_\alpha$.

\begin{lemma}[local cumulative spectral mass]\label{lem:mass}
For the tetrahedron and any measurable boundary set $F$,
\begin{equation}\label{eq:mass}
\Tr\big(\one_F\one_{[0,L]}(\Lambda_a)\one_F\big)\le C(1+L)^2|F| .
\end{equation}
\end{lemma}

\begin{proof}
The trace Sobolev inequality, the harmonic extension, and the boundary Poincar\'e inequality give
\begin{equation}\label{eq:nash}
\|f\|_4^2\le C\big(\fa_{T_a}[f]+\|f\|_2^2\big),\qquad \|f\|_2^3\le C\|f\|_1\big(\fa_{T_a}[f]+\|f\|_2^2\big).
\end{equation}
The form is Markovian: applying a normal contraction to an extension cannot increase its energy, and taking the infimum preserves this property. Apply the Nash differential inequality to the semigroup of $\Lambda_a+1$. Its $L^1$ to $L^2$ norm is at most $Ct^{-1}$; duality and composition give an $L^1$ to $L^\infty$ bound $Ct^{-2}$ for $0<t\le1$. Thus its heat diagonal is at most $Ct^{-2}$ almost everywhere. The spectral inequality $\one_{[0,L]}(\Lambda_a)\le e\,e^{-\Lambda_a/L}$ for $L\ge1$, followed by integration over $F$, proves \eqref{eq:mass}. Bounded $L$ is absorbed into the constant. This is a weighted trace statement, not a comparison of spectra after cutting the domain.
\end{proof}

The same argument on a homogeneous convex cone has no $L^2$ term in \eqref{eq:nash}. Indeed, a homogeneous bilipschitz flattening and the homogeneous half-space trace inequality give $\|f\|_4^2\le C\fa_C[f]$. Its heat bound is $Ct^{-2}$ for every $t>0$, and its spectral diagonal is bounded by $CL^2$. Both infinite dihedral wedges and trihedral tangent cones of the tetrahedron satisfy these assumptions.

For a vertex neighborhood of radius $r$, \eqref{eq:mass} is $O(L^2r^2)$. In particular a core with $r=L^{-\rho}$, $1/2<\rho<1$, has $o(L)$ total mass. This alone does not control the influence of vertex scattering on eigenfunctions outside that core.

\section{The inverse sector operator and its trace}

\subsection{The separated Neumann problem}
Let $W_\alpha=\{(r,\theta):r>0,\ 0<\theta<\alpha\}$ and $\Hcal=L^2(\R_+,dr)\oplus L^2(\R_+,dr)$. Set $k_v(r)=K_{iv}(r)$, the real Macdonald function. The Kontorovich--Lebedev transform \cite{DLMF} has normalization
\begin{equation}\label{eq:KL}
\widehat f(v)=\int_0^\infty k_v(r)f(r)\,\frac{dr}{r},\qquad f(r)=\frac{2}{\pi^2}\int_0^\infty v\sinh(\pi v)k_v(r)\widehat f(v)\,dv .
\end{equation}
It diagonalizes $-(r\partial_r)^2+r^2$ with eigenvalue $v^2$. Multiplying the massive Green equation by $r^2$ reduces the angular problem to the Neumann resolvent of $-\partial_\theta^2+v^2$ on $[0,\alpha]$. Its endpoint matrix is
\begin{equation}\label{eq:endpoint}
\frac1v\begin{pmatrix}\coth(\alpha v)&\csch(\alpha v)\\ \csch(\alpha v)&\coth(\alpha v)\end{pmatrix}.
\end{equation}
Consequently, with respect to boundary measure $dr'$,
\begin{equation}\label{eq:Bkernel}
B_\alpha(r,r')=\frac{2}{\pi^2}\int_0^\infty\sinh(\pi v)k_v(r)k_v(r')\begin{pmatrix}\coth(\alpha v)&\csch(\alpha v)\\ \csch(\alpha v)&\coth(\alpha v)\end{pmatrix}dv .
\end{equation}
Individual kernels are first understood on compactly supported smooth data, or distributionally. No unregularized diagonal is subtracted at this stage. The boundary representation of the Neumann problem has no extra factor two. At the flat angle, \eqref{eq:Bkernel} gives the kernel $\pi^{-1}K_0(|s-s'|)$ of $(1-\partial_s^2)^{-1/2}$.

\subsection{Trace-class subtraction}
The required norm identity is
\begin{equation}\label{eq:normid}
\int_0^\infty k_v(r)^2\,dr=\frac{\pi^2}{4\cosh(\pi v)} .
\end{equation}
For completeness, multiply Schl\"afli representations of $K_\nu$ and $K_{-\nu}$ and integrate first in $r$. The resulting integral is
\begin{equation}\label{eq:schlafli}
\frac{\sqrt\pi}{4}\int_0^\infty\!\!\int_0^\infty e^{-x-y}x^{-\nu-1/2}y^{\nu-1/2}(x+y)^{-1/2}\,dx\,dy .
\end{equation}
The substitution $x=qz$, $y=q(1-z)$ evaluates it as $\pi\Gamma(1/2-\nu)\Gamma(1/2+\nu)/4$. For $\nu=iv$, the gamma reflection identity gives \eqref{eq:normid}. The integrations are absolutely convergent on this imaginary axis.

The unitary ray-exchange decomposition has angular multipliers
\begin{equation}\label{eq:parity}
m_{\alpha,+}(v)=\coth(\alpha v/2),\qquad m_{\alpha,-}(v)=\tanh(\alpha v/2).
\end{equation}
In each parity sector, the difference is the rank-one integral
\begin{equation}\label{eq:rankone}
B_{\alpha,\pm}-B_{\pi,\pm}=\frac{2}{\pi^2}\int_0^\infty\sinh(\pi v)\big(m_{\alpha,\pm}(v)-m_{\pi,\pm}(v)\big)|k_v\rangle\langle k_v|\,dv .
\end{equation}
It converges in trace norm. By \eqref{eq:normid}, the integral of its rank-one trace norms is
\begin{equation}\label{eq:tracenormint}
\frac12\int_0^\infty\tanh(\pi v)\,|m_{\alpha,\pm}(v)-m_{\pi,\pm}(v)|\,dv<\infty .
\end{equation}
At zero the even difference is $O(v^{-1})$, canceled by $\tanh(\pi v)=O(v)$; the odd difference is $O(v)$. At infinity the differences decrease exponentially. The even perturbation is positive and the odd perturbation negative. Taking traces and summing yields
\begin{equation}\label{eq:Ialpha}
\begin{aligned}
I_\alpha:=\Tr(B_\alpha-B_\pi)&=\int_0^\infty\tanh(\pi v)\big(\coth(\alpha v)-\coth(\pi v)\big)\,dv\\
&=\int_0^\infty\frac{\sinh((\pi-\alpha)v)}{\cosh(\pi v)\sinh(\alpha v)}\,dv .
\end{aligned}
\end{equation}
The integrand tends to $(\pi-\alpha)/\alpha$ at zero and decreases exponentially at infinity. This proves the trace assertion of Theorem~\ref{thm:invariant} without an assertion about traces of sharp projector differences.

\subsection{Spectral shift and edge momentum}
Use the trace-class spectral-shift theorem in the convention
\begin{equation}\label{eq:ssf}
\Tr\big(f(B_\alpha)-f(B_\pi)\big)=\int f'(b)\xi_\alpha(b)\,db,\qquad \int\xi_\alpha(b)\,db=I_\alpha .
\end{equation}
Here $\xi_\alpha\in L^1$, supported in the convex hull of the spectra, and $\|\xi_\alpha\|_1\le\|B_\alpha-B_\pi\|_1$. Put $\Gamma_\alpha(\mu)=\xi_\alpha(1/\mu)$ for $\mu>0$. Then
\begin{equation}\label{eq:Gamma}
\int_0^\infty\frac{|\Gamma_\alpha(\mu)|}{\mu^2}\,d\mu<\infty,\qquad I_\alpha=\int_0^\infty\frac{\Gamma_\alpha(\mu)}{\mu^2}\,d\mu .
\end{equation}
The lower limit is zero, not one. A convex wedge has subthreshold spectrum, which must not be omitted. The flat massive cumulative projector has diagonal $\pi^{-1}\sqrt{\mu^2-1}_{\,+}$, not $\pi^{-1}(\mu-1)$.

Let $\Dcal_\alpha$ be the homogeneous DtN operator on $W_\alpha\times\R_z$. Fourier transformation in $z$ uses measure $dk/(2\pi)$. Its fiber is unitarily equivalent to $|k|D_\alpha$. For $\varphi\in C_c^\infty((0,\infty))$, apply \eqref{eq:ssf} to $f_k(b)=\varphi(k/b)$. Fubini is justified by
\begin{equation}\label{eq:fubini}
\iint_0^\infty|f_k'(b)\xi_\alpha(b)|\,dk\,db\le\|\xi_\alpha\|_1\int_0^\infty u|\varphi'(u)|\,du .
\end{equation}
Integrating in $k$ first gives the exact identity
\begin{equation}\label{eq:density}
\frac1{2\pi}\int_\R\Tr\big(\varphi(|k|D_\alpha)-\varphi(|k|D_\pi)\big)\,dk=c_\alpha\int_0^\infty\varphi(u)\,du .
\end{equation}
This already includes the continuum. A count of only the known tetrahedral guided branch, whose dispersion is $|k|/\sqrt3$, would not give \eqref{eq:density}. The next section supplies the stronger trace-norm convergence needed for spatial localization.

\section{Weighted trace ideals and ordinary local heat traces}

\subsection{A weighted inverse estimate}
Write $D_0=D_\pi$ and $B_0=B_\pi$ in this section only. Angular stretching of the half-plane by $\theta\mapsto\alpha\theta/\pi$, with $r$ fixed, preserves boundary arclength. The transformed energy has coefficients $\alpha/\pi$ on the radial and mass terms and $\pi/\alpha$ on the angular term. Minimizing over extensions gives
\begin{equation}\label{eq:stretch}
\frac\alpha\pi D_0\le D_\alpha\le\frac\pi\alpha D_0,\qquad \frac\alpha\pi B_0\le B_\alpha\le\frac\pi\alpha B_0 .
\end{equation}
Let $B_{\alpha,+}-B_{0,+}=T_+$ and $B_{\alpha,-}-B_{0,-}=-T_-$, with $T_\pm\ge0$. In the unitary cosine and sine transforms, respectively, $B_0$ is multiplication by $b(p)=(1+p^2)^{-1/2}$, and
\begin{equation}\label{eq:Fpm}
\begin{aligned}
F_+(p,v)&=\sqrt{\pi/2}\,\frac{b(p)\cos(v\arsinh p)}{\cosh(\pi v/2)},\\
F_-(p,v)&=\sqrt{\pi/2}\,\frac{b(p)\sin(v\arsinh p)}{\sinh(\pi v/2)} .
\end{aligned}
\end{equation}
These are the transforms of $k_v$. They follow by integrating its exponential representation. The apparent denominator in the odd formula has its usual continuous limiting interpretation at $v=0$.

Inserting \eqref{eq:Fpm} into the positive rank-one integrals shows that the spectral diagonal of $T_\pm$ is at most $C_\alpha b(p)^2$. Indeed, the even weight is $O(1)$ at zero, the odd weight is $O(v^2)$, and the hyperbolic quotients decrease exponentially at infinity. Consequently, for $0<q<1/2$,
\begin{equation}\label{eq:weighted}
\Tr\big(D_0^qT_\pm D_0^q\big)\le C_\alpha\int_0^\infty(1+p^2)^{q-1}\,dp<\infty .
\end{equation}
Loewner--Heinz applied to \eqref{eq:stretch} makes $D_\alpha^qD_0^{-q}$ bounded. Therefore
\begin{equation}\label{eq:Wq}
W_q:=D_\alpha^q(B_\alpha-B_0)D_0^q\in\Sone,\qquad 0<q<1/2 .
\end{equation}

\subsection{Heat-fiber integrability}
\begin{lemma}\label{lem:heatfiber}
For every $0<q<1/2$, there exist $C_q,c>0$ such that
\begin{equation}\label{eq:heatfiber}
\|e^{-uD_\alpha}-e^{-uD_0}\|_1\le C_qu^{2q-1}e^{-cu},\qquad u>0 .
\end{equation}
\end{lemma}

\begin{proof}
Put $A=D_\alpha$, $C=D_0$. The inverse-difference Duhamel identity reads
\begin{equation}\label{eq:duhamelinv}
e^{-uA}-e^{-uC}=\int_0^uAe^{-(u-v)A}(A^{-1}-C^{-1})Ce^{-vC}\,dv .
\end{equation}
Factor the integrand as $A^{1-q}e^{-(u-v)A}W_qC^{1-q}e^{-vC}$. Both generators have a positive spectral lower bound by \eqref{eq:stretch}. Functional calculus bounds the two outside factors by $C(u-v)^{q-1}e^{-c(u-v)}$ and $Cv^{q-1}e^{-cv}$. Their product is integrable in $v$, with integral $u^{2q-1}\mathrm{B}(q,q)$. For a domain justification, first replace $A,C$ by $(A^{-1}+\varepsilon)^{-1}$ and $(C^{-1}+\varepsilon)^{-1}$. The inverse difference is unchanged, the weighted ideal bound is uniform, and the integrable majorant permits passage to the limit in trace norm.
\end{proof}

In particular the trace norm in \eqref{eq:heatfiber} is integrable over $u\in(0,\infty)$. After edge Fourier transformation and scaling,
\begin{equation}\label{eq:Malpha}
\frac1{2\pi}\int_\R\|(e^{-t\Dcal_\alpha}-e^{-t\Dcal_\pi})_k\|_1\,dk=\frac{M_\alpha}{t}<\infty .
\end{equation}
Here $M_\alpha=\pi^{-1}\int_0^\infty\|e^{-uD_\alpha}-e^{-uD_0}\|_1\,du$. Its value need not equal the signed invariant $c_\alpha$.

\subsection{Spatial localization}
If a self-adjoint direct-integral operator $Q$ is translation invariant in $z$ and $\int\|Q_k\|_1\,dk<\infty$, then
\begin{equation}\label{eq:localize}
\|M_\eta QM_\eta\|_1\le\frac{\|\eta\|_2^2}{2\pi}\int\|Q_k\|_1\,dk,\qquad \eta\in L^2(\R_z)\cap L^\infty .
\end{equation}
Indeed, $|Q|^{1/2}M_\eta$ is Hilbert--Schmidt with the displayed squared norm. Factor through this operator and $\sgn Q$. Bounded transverse cutoffs can be inserted, and cutoffs tending strongly to one converge in trace norm after this localization.

A trace-class operator has an integrable diagonal in the sense of its eigenfunction expansion; the $L^1$ norm of that diagonal is bounded by the trace norm. Thus \eqref{eq:Malpha} defines a transverse density $q_\alpha\in L^1(\R)$ at time one. Homogeneity and \eqref{eq:density}, with heat testing, give
\begin{equation}\label{eq:qalpha}
Q_t(s,z;s,z)=t^{-2}q_\alpha(s/t),\qquad \int_\R q_\alpha(s)\,ds=c_\alpha,\qquad Q_t=e^{-t\Dcal_\alpha}-e^{-t\Dcal_\pi} .
\end{equation}
For a compactly supported bounded cutoff $\chi(s,z)$ continuous across the unfolded edge, dominated convergence now proves
\begin{equation}\label{eq:localheat}
\Tr\big(\chi e^{-t\Dcal_\alpha}\chi\big)=\frac{1}{2\pi t^2}\int\chi^2\,dS+\frac{c_\alpha}{t}\int\chi(0,z)^2\,dz+o(t^{-1}) .
\end{equation}
More quantitatively, if the longitudinal projection of its support has length at most $\ell$, the difference from the area term is at most $M_\alpha\ell\|\chi\|_\infty^2/t$. This establishes the ordinary localized heat assertion in Theorem~\ref{thm:invariant}.

\subsection{Smooth compact spectral tests}
The same spatial conclusion is available for a compact positive-energy test. Here is a useful uniform estimate, including its low-inverse-energy issue. In one parity sector write $B=B_0\pm T$ and $B_\theta=B_0\pm\theta T$. By \eqref{eq:stretch}, $B_\theta\ge cB_0$ uniformly. For the finite positive measure $\mu_\theta(F)=\Tr(TE_{B_\theta}(F))$, one has
\begin{equation}\label{eq:mutheta}
\mu_\theta([0,k])\le2k\Tr\big(T(B_\theta+k)^{-1}\big)\le Ck\log(2/k),\qquad 0<k<1/2 .
\end{equation}
The last inequality follows from the diagonal bound preceding \eqref{eq:weighted} and the substitution $p=\sinh v$, which reduces its integral to $\int_0^\infty(1+k\cosh v)^{-1}dv$. For $z=x+iy$, $|x|\le C_0k$, $0<|y|\le C_1k$, a dyadic division at energy $k$ gives
\begin{equation}\label{eq:resolventbound}
\Tr\big(T|B_\theta-z|^{-2}\big)\le\frac{Ck\log(2/k)}{y^2} .
\end{equation}
Factor the resolvent identity through $T^{1/2}$ and apply the Hilbert--Schmidt inequality. The trace norm of the resolvent difference has the same upper bound. For $f_k(b)=\varphi(k/b)=F(b/k)$, an almost-analytic extension in the Helffer--Sj\"ostrand formula has derivative bounded by $C_Nk^{-1}(|y|/k)^N$ on a region of area $O(k^2)$. Hence
\begin{equation}\label{eq:smoothtest}
\|\varphi(kD_\alpha)-\varphi(kD_0)\|_1\le C_\varphi\log(2/k),\qquad 0<k<1/2 .
\end{equation}
For large $k$ both operators vanish, and intermediate $k$ are handled by ordinary trace-class functional calculus. The trace norms are therefore integrable in $k$. Combining with \eqref{eq:density}, \eqref{eq:localize}, and scaling gives the localized model formula
\begin{equation}\label{eq:localtest}
\begin{aligned}
\Tr\big(\chi\varphi(h\Dcal_\alpha)\chi\big)={}&\frac{h^{-2}}{2\pi}\Big(\int_0^\infty u\varphi(u)\,du\Big)\int\chi^2\,dS\\
&+h^{-1}c_\alpha\Big(\int_0^\infty\varphi(u)\,du\Big)\int\chi(0,z)^2\,dz+o(h^{-1}) .
\end{aligned}
\end{equation}
The constants in \eqref{eq:smoothtest} involve derivatives of $\varphi$. This estimate must not be asserted uniform for a spectral test whose transition width shrinks with $h$.

\section{Spectral-window estimates on the compact tetrahedron}

\subsection{Exact cone window bounds from scaling}
Let $B=\Lambda_C$ on an infinite dihedral cone or tetrahedral tangent cone. Its spectral diagonal satisfies $0\le e_B(y,y;L)\le CL^2$ by Section~2. In the dihedral case, translation invariance and homogeneity imply
\begin{equation}\label{eq:scaling}
e_B((s,z),(s,z);L)=L^2F(Ls),\qquad 0\le F\le C .
\end{equation}
Put $H(R)=\int_{-R}^RF(u)\,du$. The cumulative strip trace per unit edge length is $a_L(\gamma)=LH(L\gamma)$. For $L_2>L_1$,
\begin{equation}\label{eq:stripincrement}
\begin{aligned}
0\le a_{L_2}(\gamma)-a_{L_1}(\gamma)&=(L_2-L_1)H(L_2\gamma)+L_1\big(H(L_2\gamma)-H(L_1\gamma)\big)\\
&\le C\gamma(L_2^2-L_1^2) .
\end{aligned}
\end{equation}
In a trihedral cone, the trace in a radius-$r$ ball equals the integral of $e_B(u,u;1)$ over the boundary ball of radius $Lr$. Its increment is bounded by $Cr^2(L_2^2-L_1^2)$. Positivity of interval projectors then implies, for a contraction $M$ supported in the respective set,
\begin{equation}\label{eq:dtrace}
d\Tr\big(M^*\one_{[0,t]}(B)M\big)\le CVt\,dt,\qquad V=|I|\gamma\ \text{or}\ r^2 .
\end{equation}
These interval inequalities do not differentiate a pointwise spectral kernel. Homogeneity and the bounded positive diagonal suffice, so all model channels are included.

\subsection{A bounded intertwiner}
Fix the infinite wedge or trihedral cone $C$ containing $T_a$ at the edge or vertex under consideration. For the $L^2$ Dirichlet solution on this global Lipschitz graph domain, define
\begin{equation}\label{eq:intertwiner}
Rf=(H_Cf)|_{\partial T_a} .
\end{equation}
This is bounded from $L^2(\partial C)$ to $L^2(\partial T_a)$. On common faces it is the identity. On each remaining plane, surface measure inside $C$ is a Carleson measure: a boundary tent of radius $r$ contains at most $Cr^2$ of that surface measure. The Carleson embedding inequality and the classical $L^2$ Dirichlet nontangential-maximal estimate give
\begin{equation}\label{eq:carleson}
\int_{\partial T_a\setminus\partial C}|H_Cf|^2\,dS\le C\|N_C(H_Cf)\|_{L^2(\partial C)}^2\le C\|f\|_2^2 .
\end{equation}
Let $\eta$ be a fixed smooth ambient cutoff supported only on common faces, equal to one on a smaller patch. Write $A=\Lambda_a$, $B=\Lambda_C$, and $J=\eta R$. On its support, $J$ is simply multiplication followed by identification of the common boundary. Uniqueness of harmonic extension gives $\eta ARf=JBf$, initially on regular data. Therefore
\begin{equation}\label{eq:intertwining}
AJ-JB=[A,\eta]R=:K .
\end{equation}
The $L^2$ DtN commutator theorem on a Lipschitz domain \cite[Theorem~1.1]{S} makes $[A,\eta]$ bounded. Thus $K$ is bounded. The identity extends weakly and then by density; the classical $L^2$ regularity problem supplies the operator-domain interpretation. The cutoff $\eta$ is fixed before the smaller strip or ball is chosen. This fact is essential to uniformity.

\subsection{An abstract transfer lemma}
\begin{lemma}\label{lem:transfer}
Let $A\ge0$ have compact resolvent with $\rank\one_{[0,L]}(A)\le C(1+L)^2$. Let $B\ge0$ be self-adjoint, and suppose $J$ and $K=AJ-JB$ are bounded. If the contraction $M$ satisfies \eqref{eq:dtrace}, then for fixed $w>0$,
\begin{equation}\label{eq:transfer}
\|\one_{[L-w,L+w]}(A)JM\|_{\Stwo}^2\le C(VL+1) .
\end{equation}
For $1\le w\le L/4$, the right-hand side is $C(wVL+1)$.
\end{lemma}

\begin{proof}
Put $E=\one_{[L-w,L+w]}(A)$ and $Q=\one_{[L/2,2L]}(B)$. For $S=(B-L-i)^{-1}Q$, the intertwining identity implies
\begin{equation}\label{eq:EJQM}
EJQM=E(A-L-i)JSM-EKSM .
\end{equation}
By \eqref{eq:dtrace},
\begin{equation}\label{eq:SM}
\|SM\|_2^2\le CV\int_{L/2}^{2L}\frac{t}{(t-L)^2+1}\,dt\le CVL .
\end{equation}
The outside factor in \eqref{eq:EJQM} is bounded by $\sqrt{w^2+1}\,\|J\|+\|K\|$. For the off-band part, set $S_o=(B-L)^{-1}(1-Q)$ and $X=EJ(1-Q)$. Since $\|S_o\|\le2/L$,
\begin{equation}\label{eq:offband}
X=(A-L)XS_o-EKS_o,\qquad \|X\|\le4\|K\|/L
\end{equation}
when $L>4w$. The range of $X$ lies in the at most $C(1+L+w)^2$ dimensional range of $E$. Hence $\|XM\|_2^2\le C$. This avoids the logarithmic divergence that would result from treating the \emph{full} spatially localized resolvent as Hilbert--Schmidt in boundary dimension two. Summing band and off-band bounds proves \eqref{eq:transfer}. For wide windows replace $i$ by $iw$ in $S$; the integral in \eqref{eq:SM} is then $O(VL/w)$ and the squared prefactor $O(w^2)$. The off-band estimate is unchanged.
\end{proof}

On the smaller common patch where $\eta=1$, the actual multiplier $\widetilde M$ and model multiplier $M$ satisfy $JMM^*J^*=\widetilde M\widetilde M^*$. Therefore the actual weighted spectral trace equals $\|EJM\|_2^2$. Apply Lemma~\ref{lem:transfer} with \eqref{eq:dtrace} and the cumulative rank estimate from Section~2. This proves Theorem~\ref{thm:windows}.

\begin{remark}[Limitation]
A bounded intertwiner transfers an upper bound, not the linear coefficient. Replacing a first-order operator $B$ in boundary dimension two by $B+c$ preserves a bounded intertwining defect but changes $AL^2$ to $AL^2-2AcL+O(1)$. Thus \eqref{eq:intertwining} must not be promoted to a coefficient comparison.
\end{remark}

\section{Quantitative departure of the wedge continuum}

This section concerns the exact infinite wedge. It does not assert a finite-tetrahedron propagation theorem.

\subsection{Dilation regularity of the inverse perturbation}
Write $D=D_\alpha$, $B=D^{-1}$, $B_0=(1-\partial_s^2)^{-1/2}$, and $K=B-B_0$. Let $Q=-i(s\partial_s+1/2)$ be the self-adjoint dilation generator on the joined boundary line. Its unitary group is $U_bf(s)=b^{1/2}f(bs)$.

\begin{lemma}\label{lem:dilation}
For all nonnegative integers $n,m$, both $Q^nKQ^m$ and $|s|^nK|s|^m$ are trace class. The same holds for finite mixed words in $|s|$ and $Q$.
\end{lemma}

\begin{proof}
The Bessel square moments satisfy
\begin{equation}\label{eq:besselmoments}
\frac{\int_0^\infty r^{2n}k_v(r)^2\,dr}{\|k_v\|_2^2}=\frac{\binom{2n}{n}}{4^n}\prod_{j=0}^{n-1}\big(v^2+(j+1/2)^2\big) .
\end{equation}
This follows from the gamma integral for two Macdonald functions, or its moment recurrence. The differential equation and integration by parts give
\begin{equation}\label{eq:besselibp}
\|r^p(r\partial_r)k_v\|_2^2+\|r^{p+1}k_v\|_2^2=\Big(v^2+\frac{(2p+1)^2}{2}\Big)\|r^pk_v\|_2^2 .
\end{equation}
The differential equation reduces higher dilation derivatives to these quantities. Each fixed word therefore costs only a polynomial in $v$. The trace-norm majorant \eqref{eq:tracenormint} is integrable with every polynomial weight, proving the assertion by the same rank-one integral.
\end{proof}

Both $B$ and $B_0$ belong to $C^\infty(Q)$. Resolvent commutators and an almost-analytic functional calculus then give
\begin{equation}\label{eq:Cinfty}
Q^n\big(f(B)-f(B_0)\big)Q^m\in\Sone,\qquad f\in C_c^\infty(\R) .
\end{equation}
Indeed, insert the resolvent identity, commute the outside powers through the resolvents, and use Lemma~\ref{lem:dilation} on the middle factor. Each resulting trace norm has a finite polynomial singularity in the imaginary part of the spectral parameter, absorbed by the order of the almost-analytic extension. The corresponding statement holds for positive compactly supported functions of $D$ by substituting $b\mapsto1/b$.

\subsection{Mourre estimate and possible embedded point spectrum}
The free identities are
\begin{equation}\label{eq:freecomm}
i[B_0,-Q]=B_0-B_0^3,\qquad i[D_0,Q]=D_0-D_0^{-1} .
\end{equation}
Lemma~\ref{lem:dilation} gives $i[B,-Q]=B-B^3+T$ with $T$ trace class. On a compact interval $I\Subset(1,\infty)$,
\begin{equation}\label{eq:mourre}
E_I(D)i[D,Q]E_I(D)=E_I(D)(D-D^{-1})E_I(D)+\big(DE_I(D)\big)T\big(E_I(D)D\big) .
\end{equation}
The localized commutator has a bounded form interpretation; $DTD$ is not asserted globally bounded. The compact-error Mourre theorem, with arbitrarily many cutoff commutators, gives absence of singular continuous spectrum in $(1,\infty)$ and locally finite point spectrum of finite multiplicity there.

One can remove this point spectrum without assuming its absence. If $Bu=bu$, $0<b<1$, put $F=-Ku$. Lemma~\ref{lem:dilation} gives every spatial moment and dilation derivative of $F$. In Fourier variables,
\begin{equation}\label{eq:eigeneq}
\big((1+p^2)^{-1/2}-b\big)\widehat u(p)=\widehat F(p) .
\end{equation}
The multiplier has two simple real zeros. Since $\widehat u\in L^2$ and $\widehat F$ is smooth, $\widehat F$ vanishes at both. Locally, division by a simple zero preserves all finite Sobolev orders with the loss of one derivative, as follows from $g(x)/x=\int_0^1g'(tx)\,dt$ when $g(0)=0$. Away from the zeros, the reciprocal and all its derivatives are bounded. Hence $|s|^Nu\in L^2$ for every $N$; the same argument with $p\partial_p$ gives $u\in D(Q^N)$ for every $N$.

Choose a larger compact interval containing the desired energy support, and let $P$ be the sum of its eigenspace projections. It is finite rank, spatially rapidly decreasing, and dilation regularizing. The operator $Q_c=(1-P)Q(1-P)$ is self-adjoint on $D(Q)$ and differs from $Q$ by a bounded regularizing finite-rank operator. On the complementary subspace, compression of the compact error in \eqref{eq:mourre} to sufficiently small energy intervals tends to zero in norm. Thus a strict Mourre estimate holds on each member of a finite energy cover.

\subsection{From dilation escape to spatial escape}
Let $X_R$ be any multiplication contraction supported in $|s|\le R$. For $\eta\in C_c^\infty((1,\infty))$,
\begin{equation}\label{eq:QNeta}
\|Q^N\eta(D)X_R\|\le C_{\eta,N}(1+R)^N .
\end{equation}
To see this, subtract the free cutoff using \eqref{eq:Cinfty}. The free cutoff is convolution by a Schwartz kernel $k$. Acting with $Q^N$ produces sums of $s^jk^{(j)}(s-s')X_R(s')$. The inequality $|s|^j\le C_j(R^j+|s-s'|^j)$ and Schur's test prove \eqref{eq:QNeta}. The same estimate holds with $Q_c$ and $(1-P)$ inserted.

We use the minimal-escape theorem of Hunziker--Sigal--Soffer \cite[Theorem~1.1]{HSS} in its precise form: a strict Mourre estimate $E_\Delta(H)i[H,A]E_\Delta(H)\ge\theta E_\Delta(H)$ and bounded cutoff commutators of every order imply
\begin{equation}\label{eq:HSS}
\big\|\one_{(-\infty,b+\vartheta t]}(A)e^{-itH}g(H)\one_{[b,\infty)}(A)\big\|\le C_Nt^{-N},\qquad 0<\vartheta<\theta,
\end{equation}
uniformly in $b\in\R$, for $g\in C_c^\infty(\Delta)$ and $t\ge1$. All its hypotheses have been verified on the complementary space. Insert the dilation projection $\one_{[-L,L]}(Q_c)$ on both sides of the energy-localized propagator, with $L=\vartheta t/4$. By \eqref{eq:QNeta}, each discarded tail costs $C_N((1+R)/t)^N$. For the middle term, apply \eqref{eq:HSS} with $b=-L$. A finite partition over the strict Mourre intervals, and adjoints for negative time, prove
\begin{equation}\label{eq:spatialescape}
\|X_Re^{itD}g(D)P_{\mathrm{ac}}X_R\|\le C_{g,N}\Big(\frac{1+R}{1+|t|}\Big)^N,\qquad g\in C_c^\infty((1,\infty)) .
\end{equation}
This is a spatial operator-norm estimate, stronger than qualitative scattering completeness or time-averaged RAGE decay.

For a shrinking energy window $\varphi((D-\tau)/\rho)$ contained in a fixed compact subinterval of $(1,\infty)$, Fourier convolution of \eqref{eq:spatialescape} gives
\begin{equation}\label{eq:shrinkingwindow}
\begin{aligned}
\|X_Re^{itD}\varphi((D-\tau)/\rho)P_{\mathrm{ac}}X_R\|\le{}&C_N\Big(\frac{1+R}{1+|t|}\Big)^N\\
&+C_M(1+\rho|t|)^{-M} .
\end{aligned}
\end{equation}
The proof splits the convolution at $|u|=|t|/2$. In the inner part use \eqref{eq:spatialescape}; in the tail use rapid decay of $\rho\widehat\varphi(\rho u)$. No uncontrolled derivatives of a shrinking window enter the constants.

\subsection{Semiclassical product-wedge consequence}
For mass $m>0$, exact dilation gives $D_m=mU_mDU_m^{-1}$. In edge Fourier variables put $k=h\zeta$, and restrict to a fixed annulus $\kappa\le|k|\le1-\epsilon$. If $\psi_\gamma$ is supported at distance at most $C\gamma$ from the edge and $\Pi_c$ is the direct integral of the fiber absolutely continuous projections, \eqref{eq:shrinkingwindow} gives
\begin{equation}\label{eq:semiclassical}
\begin{aligned}
&\big\|a(hD_z)\psi_\gamma e^{iT\Dcal_\alpha}\varphi((h\Dcal_\alpha-1)/\eta)\Pi_c\psi_\gamma\big\|\\
&\qquad\le C_N\Big(\frac{\gamma+h/\kappa}{T}\Big)^N+C_M\Big(\frac{h}{\eta T}\Big)^M .
\end{aligned}
\end{equation}
The model parameters are $t=|k|T/h$, $R\le C|k|\gamma/h$, $\tau=1/|k|$, and $\rho=\eta/|k|$. They remain in fixed compact continuum bands for fixed $\kappa,\epsilon$ and sufficiently small $\eta$. This proves \eqref{eq:semiclassical} uniformly in the displayed semiclassical scales. It is not an endpoint-uniform estimate as $\kappa\downarrow0$ or $\epsilon\downarrow0$, and it is not itself a locality theorem for the finite tetrahedron.

\section{Static locality and short-time cumulative traces}

This section distinguishes three assertions which are not interchangeable: a bounded intertwining defect, a smooth spectral trace comparison, and a cumulative trace comparison at a fixed physical time. The first was proved in Section~5. The second requires an operator-ideal improvement, and the third additionally uses upper propagation. None uses a lower escape estimate at zero tangential momentum.

\subsection{Hilbert--Schmidt defects from the common harmonic problem}
Retain $A=\Lambda_a$, the half-space or product-wedge model $B$, and the common-patch identification $J=\eta R$ of Section~5. Choose nested real cutoffs $\chi\prec\zeta\prec\eta\prec\eta_1$ in a fixed vertex-free common patch; the notation means that each later cutoff equals one near the preceding support. The inner multiplier $\chi$ may be any multiplication contraction supported in a fixed smaller patch. Enclosing cutoffs, not derivatives of $\chi$, determine the constants. Put $K=AJ-JB$.

\begin{lemma}\label{lem:HS}
The following operators are Hilbert--Schmidt:
\begin{equation}\label{eq:HSops}
\rho D(1-\widetilde\rho),\qquad K\zeta,\qquad (1-\eta_1)K,\qquad D=A\text{ or }B,
\end{equation}
where $\rho\prec\widetilde\rho$ are compact common-patch cutoffs. Consequently,
\begin{equation}\label{eq:HSK}
\|\one_{[0,L]}(A)K\|_{\Stwo}\le C(1+L) .
\end{equation}
\end{lemma}

\begin{proof}
For remote boundary data $f=(1-\widetilde\rho)g$, the harmonic extension $u$ has zero Dirichlet data in an enclosing patch. The $L^2$ Dirichlet estimate, local Caccioppoli inequality and mean-value inequality give $\sup|u|\le C\|g\|_2$ on a smaller patch. One can apply the mean-value inequality up to the zero-data boundary by extending $|u|^2$ by zero; this extension is locally subharmonic. On an open face, odd reflection bounds $|\partial_\nu u|$ by $C\|g\|_2$.

At an edge let $\nu=\pi/\alpha_0>2$, and choose $1/2<\beta<\gamma<\nu$. In normalized coordinates $|z|<1$, put $d(z)=1-z^2$. In the cusp $0<r<\varepsilon d(z)$ the function
\begin{equation}\label{eq:barrier}
w(r,\theta,z)=r^\beta d(z)^{-\beta}\cos\big(\gamma(\theta-\alpha_0/2)\big)
\end{equation}
is positive and superharmonic for sufficiently small fixed $\varepsilon$. Indeed its transverse Laplacian is $(\beta^2-\gamma^2)r^{\beta-2}d^{-\beta}$ times the angular factor, whereas its longitudinal second derivative is bounded by $C_\beta r^\beta d^{-\beta-2}$ times that factor. On $r=\varepsilon d(z)$ it has a fixed positive lower bound. The maximum principle, first for regular remote data and then by approximation, gives $|u|\le Cr^\beta\|g\|_2$ in a smaller longitudinal interval. Odd reflection in a ball of radius comparable to $r$ centered on either face gives
\begin{equation}\label{eq:normalderiv}
|\partial_\nu u(r,z)|\le Cr^{\beta-1}\|g\|_2 .
\end{equation}
The square of this bound is integrable against $dr\,dz$. An operator satisfying $|Tf(y)|\le w(y)\|f\|_2$ with $w\in L^2$ is Hilbert--Schmidt: apply the bound to finite linear combinations of an orthonormal basis, sum the squared values pointwise, and integrate. This proves the first assertion of \eqref{eq:HSops}, including its adjoint.

For data $\zeta f$, write $R\zeta=J\zeta+S$. The operator $S$ is harmonic restriction to the nonshared faces. Those faces are separated from the support of $\zeta$, so the same pointwise harmonic bound makes $S$ Hilbert--Schmidt. Since $\eta=1$ near $\zeta$,
\begin{equation}\label{eq:Kzeta}
K\zeta=(1-\eta)AJ\zeta+[A,\eta]S,\qquad (1-\eta_1)K=(1-\eta_1)A\eta R .
\end{equation}
Separated blocks and bounded commutators prove the remaining assertions. Finally decompose $K=\eta_1K+K_\infty$, with $K_\infty\in\Stwo$, and apply the local cumulative bound \eqref{eq:mass} to $\eta_1$. This proves \eqref{eq:HSK}.
\end{proof}

The model graph estimate $D(B)=H^1(\partial C)$, with $\|f\|_{H^1}\le C(\|Bf\|_2+\|f\|_2)$, follows from the harmonic Rellich identity with a constant vector having a fixed nonzero scalar product with both face normals. The weak identity and closure supply the operator-domain statement. Bounded commutators $[B,\rho]$ follow from the corresponding Lipschitz-graph estimate, or from \eqref{eq:intertwining}, the compact-domain commutator theorem, and the separated Hilbert--Schmidt blocks just proved. These are first-order elliptic estimates; no propagation conclusion is being inferred from them alone.

\subsection{Smooth functional calculus with an $O(1)$ local trace error}
\begin{lemma}\label{lem:smoothcalc}
For $\psi\in C_c^\infty((0,\infty))$,
\begin{equation}\label{eq:psiHS}
\|(\psi(hA)J-J\psi(hB))\chi\|_{\Stwo}\le C_\psi h,\qquad 0<h\le1 .
\end{equation}
Consequently, for $\varphi\in C_c^\infty((0,\infty))$,
\begin{equation}\label{eq:smoothtrace}
\Tr\big(\chi_A\varphi(hA)\chi_A\big)-\Tr\big(\chi_B\varphi(hB)\chi_B\big)=O_\varphi(1) .
\end{equation}
Constants are controlled by finitely many smooth seminorms on a fixed energy-support set, and are uniform for the bounded inner multipliers specified above.
\end{lemma}

\begin{proof}
Write $R_h(z)=(hB-z)^{-1}$. The graph estimate and Fourier truncation in two-dimensional compact coordinate boxes give
\begin{equation}\label{eq:singvals}
s_n(\rho R_h(z))\le C|\operatorname{Im}z|^{-1}\min(1,h^{-1}n^{-1/2}),\qquad \|\rho R_h(z)\|_{\mathfrak S_4}\le Ch^{-1/2}|\operatorname{Im}z|^{-1} .
\end{equation}
A once-localized order-minus-one resolvent is not claimed to be Hilbert--Schmidt in dimension two.

For separated cutoffs $\chi\prec\theta$, almost-analytic functional calculus and two intermediate commutators give
\begin{equation}\label{eq:separated}
\|(1-\theta)f(hB)\chi\|_{\Stwo}\le C_fh .
\end{equation}
Here are the details of the two-commutator step. Choose $\chi\prec\rho_0\prec\rho_1\prec\rho_2\prec\theta$ and set $m=1-\theta$, $K_j=[B,\rho_j]$. The first commutator writes $mf(hB)\chi=m[f(hB),\rho_0]\chi$. Split
\begin{equation}\label{eq:K0split}
K_0=\rho_1K_0\rho_1+T,\qquad T=(1-\rho_1)B\rho_0-\rho_0B(1-\rho_1)\in\Stwo .
\end{equation}
In the functional-calculus integral $hmR_hK_0R_h\chi$, the $T$ term is $O_{\Stwo}(h|\operatorname{Im}z|^{-2})$. For the localized term insert $mR_h\rho_1=-hmR_hK_2R_h\rho_1$. The product of its two compactly localized resolvents is in $\Stwo$ by \eqref{eq:singvals}, with norm $O(h^{-1}|\operatorname{Im}z|^{-2})$. The two factors of $h$ leave $O(h|\operatorname{Im}z|^{-3})$. An almost-analytic extension of sufficiently high order makes both estimates integrable. This proves \eqref{eq:separated}.

To prove \eqref{eq:psiHS}, Fourier functional calculus writes $\psi(hA)J-J\psi(hB)$ as a linear function $\mathcal L_{h,\psi}(K)$. On Hilbert--Schmidt inputs its norm is at most $C_\psi h\|K\|_2$. Split $K=K\zeta+K(1-\zeta)$. Lemma~\ref{lem:HS} handles the first part. In a normalized eigenbasis $Ae_n=a_ne_n$, a row of the remaining part is
\begin{equation}\label{eq:row}
\Big\langle K^*e_n,(1-\zeta)\frac{\psi(ha_n)-\psi(hB)}{a_n-B}\chi\,\cdot\Big\rangle .
\end{equation}
The divided difference is continuous at $B=a_n$. Because $(1-\zeta)\chi=0$, its norm is bounded by the commutator of that divided difference with $\zeta$. Fourier inversion gives a uniform bound $C_\psi h^2\|[B,\zeta]\|$: represent the divided difference of exponentials by an integral, and commute the second exponential with $\zeta$. Summing rows for $a_n\le C/h$, \eqref{eq:HSK} gives $Ch^2\|\one_{[0,C/h]}(A)K\|_2\le Ch$.

For $a_n>C/h$, with $C$ above the energy support of $\psi$, the separated-energy divided difference has the convergent representation
\begin{equation}\label{eq:highenergy}
-\int_0^\infty e^{-sA_H}E_HK(1-\zeta)e^{sB}\psi(hB)\chi\,ds .
\end{equation}
Put $s=hv$. Estimate \eqref{eq:separated} applied to $e^{v\cdot}\psi$ bounds the model factor by $Ch(1+v)^Ne^{Mv}$, whereas $e^{-hvA_H}$ has norm at most $e^{-Cv}$, with $C>M$. Integration gives $O(h^2)$. This proves \eqref{eq:psiHS}.

For $\varphi=\psi^2$, put $X=\psi(hA)J\chi$ and $Y=J\psi(hB)\chi$. The actual trace is $\|X\|_2^2$, while the model trace is $\|\psi(hB)\chi\|_2^2$. The difference is
\begin{equation}\label{eq:tracediff}
2\operatorname{Re}\langle Y,X-Y\rangle_{\Stwo}+\|X-Y\|_2^2-\|(1-\eta^2)^{1/2}\psi(hB)\chi\|_2^2 .
\end{equation}
The three terms are respectively $O(1)$, $O(h^2)$ and $O(h^2)$ by the local spectral-mass bound, \eqref{eq:psiHS}, and \eqref{eq:separated}. A general real test is a difference of two smooth squares: if $\rho=1$ near its support, use $\varphi=((\rho+\varphi)^2-(\rho-\varphi)^2)/4$. Complex tests follow by linearity. This proves \eqref{eq:smoothtrace}.
\end{proof}

The same proof uses no derivatives of the innermost multiplier. Where a localized resolvent is needed, factor it through a fixed smooth cutoff which equals one on the permitted inner support. This establishes the asserted uniformity, including for shrinking edge collars.

\subsection{Heat testing and normalized annuli}\label{subsec:heat}
A smooth dyadic partition on positive energies extends \eqref{eq:smoothtrace} to the heat function:
\begin{equation}\label{eq:heatcompare}
\Tr\big(\chi_Ae^{-tA}\chi_A\big)-\Tr\big(\chi_Be^{-tB}\chi_B\big)=O(1+|\log t|) .
\end{equation}
Indeed, on $u\asymp2^j$ write $e^{-u}\rho(2^{-j}u)=g_j(2^{-j}u)$ with $g_j(v)=\rho(v)e^{-2^jv}$ and use parameter $t/2^j$. There are $O(1+|\log t|)$ bands with $t\le2^j\le1$; their seminorms are bounded. Higher bands have exponentially summable seminorms, and bounded physical frequencies have bounded local mass. More generally, the local comparison is $O(1+|\log h|)$ for every fixed function smooth on $[0,\infty)$ and rapidly decreasing at infinity.

All these estimates are uniform on normalized conical annuli of the tetrahedron. To verify this explicitly, dilate about a vertex by $x=v+rX$. The normalized face/edge charts and enclosing separations range over a finite fixed family. Dirichlet constants are scale invariant; the factor $r$ in the normalized DtN operator cancels the $r^{-1}$ in a physical cutoff's Lipschitz norm. Local Nash bounds remain uniform. Formula \eqref{eq:HSK} uses localized spectral mass, not the growing total area of the normalized tetrahedron.

For the harmonic restriction $S$ in \eqref{eq:Kzeta}, growing nonshared faces are controlled by decay rather than their area. Compactly supported normalized model data satisfy
\begin{equation}\label{eq:decay}
|H_C(\zeta f)(x)|\le C(1+|x|)^{-2}\|f\|_2\quad\text{outside a fixed ball.}
\end{equation}
Choose one bounding half-space $n\cdot x>0$. On a fixed enclosing sphere the zero-data face estimates, with $1<\beta<\pi/\alpha_0$, bound the solution by $C(n\cdot x)\|f\|_2$. Compare in the exterior cone with the positive harmonic dipole $(n\cdot x)/|x|^3$; it is nonnegative on the other faces. The exterior maximum principle gives \eqref{eq:decay}. Its square has bounded integral over any plane in three dimensions. Hence $S$ has a uniform Hilbert--Schmidt bound even when the nonshared faces expand.

For an annular cutoff $\psi_r$, the model heat estimate from Section~4, \eqref{eq:heatcompare}, and scaling therefore give
\begin{equation}\label{eq:annulusheat}
\Big|\Tr\big(\psi_re^{-tA}\psi_r\big)-\frac1{2\pi t^2}\int\psi_r^2\,dS\Big|\le C(1+r/t),\qquad r\ge t .
\end{equation}
Use a squared dyadic partition of a radius-$\delta$ vertex neighborhood down to $r\asymp t$, and \eqref{eq:mass} on the wavelength-sized core. Area terms and traces add exactly, yielding
\begin{equation}\label{eq:vertexheat}
\Big|\Tr\big(\zeta_\delta e^{-tA}\zeta_\delta\big)-\frac1{2\pi t^2}\int\zeta_\delta^2\,dS\Big|\le C\delta/t+C\log(2+\delta/t)+C .
\end{equation}
On the fixed vertex-free complement, \eqref{eq:heatcompare} transfers the complete localized model heat formula \eqref{eq:localheat}. The omitted edge length is $O(\delta)$. First let $t\downarrow0$ with $\delta$ fixed, and then let $\delta\downarrow0$. This proves \eqref{eq:heat} without Ivrii's local edge expansion. A global logarithmic remainder is not asserted: the model-tail limit and the removal of vertices are two-limit arguments.

\subsection{From the static estimate to a fixed-time cumulative comparison}
Fix a real even $\kappa\in C_c^\infty([-1,1])$, equal to one near zero, and use the Fourier convention $\widehat\rho_S(t)=\kappa(t/S)$. Thus $\rho_S(u)=S\rho_1(Su)$, $\int\rho_S=1$, and every positive-order moment vanishes. The kernel may be signed. Set
\begin{equation}\label{eq:FS}
F_S(x)=\int_{-\infty}^x\rho_S(u)\,du=F_1(Sx),\qquad N^S_{D,\chi}(L)=\Tr\big(\chi F_S(L-D)\chi\big) .
\end{equation}
The latter is the convolution of the localized cumulative count with $\rho_S$. Quadratic local spectral mass and Schwartz tails make this an ordinary convergent trace.

\begin{proposition}[fixed-time transfer]\label{prop:fixedtime}
There is a fixed $S_0>0$, depending only on the enclosing common chart, such that
\begin{equation}\label{eq:fixedtime}
N^{S_0}_{A,\chi}(L)-N^{S_0}_{B,\chi}(L)=o(L) .
\end{equation}
The same $S_0$ and the same little-o majorant work for every multiplication contraction supported in a fixed smaller patch. In particular, $S_0$ does not decrease with an inner edge-collar width. On normalized conical annuli the comparison is uniformly $O(1+L)$; uniform little-o over that changing family is not asserted.
\end{proposition}

\begin{proof}
We use only upper, not lower, propagation. In its precise needed form, for fixed separated compact model cutoffs $\chi_0$ and $\tau$, and $\widetilde\psi\in C_c^\infty((0,\infty))$, sufficiently small fixed $S_0$ gives
\begin{equation}\label{eq:upperlocality}
\sup_{|t|\le S_0}\|\tau e^{-itB}\widetilde\psi(hB)\chi_0\|=O(h^N)\quad\text{for every }N .
\end{equation}
This is the all-angle upper-locality theorem \cite[Propositions~3.1 and~3.5]{I}, in energy-localized operator form. The initial and target supports have a positive buffer beyond the upper-speed dependence region. Fourier inversion of its time-Fourier form, with a slightly larger time cutoff and a fixed energy cutoff, gives \eqref{eq:upperlocality}; arbitrary powers absorb the Fourier normalization, and the cutoff buffer controls convolution tails. There is no restriction excluding zero tangential momentum. The ordinary smooth-face version follows from the real-principal-type parametrix. Both forms are used only with fixed enclosing cutoffs.

Choose $\chi\prec\chi_0\prec\chi_1\prec\eta\prec\chi_2$ and put $\tau=\chi_2-\chi_1$. Besides $K\chi_1\in\Stwo$, one has
\begin{equation}\label{eq:Ksplit}
K(1-\chi_2)=-JB(1-\chi_2)\in\Stwo,\qquad K=F+K\tau,\qquad F=K\chi_1+K(1-\chi_2)\in\Stwo .
\end{equation}
The extra identity uses $J(1-\chi_2)=0$ and the separated model block of Lemma~\ref{lem:HS}. It allows \eqref{eq:upperlocality} to be used only with a compact target; no assertion about the entire unbounded distant boundary is needed.

Choose real $\psi,\widetilde\psi\in C_c^\infty((0,\infty))$ with $\psi=1$ near one and $\widetilde\psi=1$ near its support. Put $Z_h=\psi(hB)\chi$. Its Hilbert--Schmidt norm is $O(h^{-1})$. Split
\begin{equation}\label{eq:Zsplit}
\tau e^{-itB}Z_h=\tau e^{-itB}\widetilde\psi(hB)\chi_0Z_h+\tau e^{-itB}\widetilde\psi(hB)(1-\chi_0)Z_h .
\end{equation}
The first term is $O_{\Stwo}(h^{N-1})$ by \eqref{eq:upperlocality}; the second is $O_{\Stwo}(h)$ by \eqref{eq:separated}. Since $\psi(hB)\to0$ strongly and $F\in\Stwo$, finite-rank approximation gives $\|F\psi(hB)\|_2\to0$. Thus
\begin{equation}\label{eq:epsh}
\sup_{|t|\le S_0}\|Ke^{-itB}Z_h\|_2\le\|F\psi(hB)\|_2+Ch=:\epsilon(h)\longrightarrow0 .
\end{equation}
For $V(t)=e^{-itA}J-Je^{-itB}$, bounded-quasicommutator Duhamel yields
\begin{equation}\label{eq:Vt}
\|V(t)Z_h\|_2\le|t|\epsilon(h),\qquad |t|\le S_0 .
\end{equation}
The integral is Hilbert--Schmidt-valued by \eqref{eq:epsh}. No decay rate for the compact spectral tail is needed.

Define $\mathcal W_D(h,t)=\Tr(\chi\psi(hD)^2e^{-itD}\chi)$. With $X_h=\psi(hA)J\chi$ and $Y_h=JZ_h$, \eqref{eq:psiHS} gives $\|X_h-Y_h\|_2=O(h)$ while all three factors have Hilbert--Schmidt norm $O(h^{-1})$. Replacing $X_h$ by $Y_h$ in the actual wave trace costs $O(1)$. The remaining difference is
\begin{equation}\label{eq:wavediff}
\Tr\big(Y_h^*V(t)Z_h\big)+\Tr\big(Z_h^*(\eta^2-I)e^{-itB}Z_h\big) .
\end{equation}
The first term is $O(h^{-1}|t|\epsilon(h))$. For the second, $\|(1-\eta^2)^{1/2}Z_h\|_2=O(h)$ by \eqref{eq:separated}, whereas the other factor has norm $O(h^{-1})$; it is $O(1)$. Hence
\begin{equation}\label{eq:wavetrace}
|\mathcal W_A(h,t)-\mathcal W_B(h,t)|\le C+Ch^{-1}|t|\epsilon(h) .
\end{equation}
Set $h=L^{-1}$. The difference between smoothing times $S_0$ and $h$, on this energy band, is the Fourier integral with multiplier
\begin{equation}\label{eq:multiplier}
\frac1{2\pi}e^{it/h}\frac{\kappa(t/S_0)-\kappa(t/h)}{it} .
\end{equation}
It vanishes for $|t|\le ch$. Integrating \eqref{eq:wavetrace} against \eqref{eq:multiplier} costs $O(1+|\log h|)$ for the constant term, but only $O(h^{-1}\epsilon(h))$ for the second term: the Duhamel factor $|t|$ cancels $1/t$. In particular no logarithm multiplies $\epsilon(h)/h$.

To account for all other energies, first note that $F_h(h^{-1}-D)=f_0(hD)$ with $f_0(u)=F_1(1-u)$. The extension of \eqref{eq:smoothtrace} following \eqref{eq:heatcompare} compares these traces with error $O(1+|\log h|)$. Next put $q=1-\psi^2$. Off the fixed energy band,
\begin{equation}\label{eq:offband2}
\begin{aligned}
&q(u)\big[F_1((S_0/h)(1-u))-F_1(1-u)\big]=g(u)+r_h(u),\\
&g(u)=q(u)\big[\one_{u<1}-F_1(1-u)\big],\qquad |r_h(u)|\le C_Mh^M(1+u)^{-M} .
\end{aligned}
\end{equation}
The apparent jump of $g$ is removed on an entire neighborhood by $q=0$. It is another fixed smooth Schwartz-tail test. Local spectral mass bounds make the trace of $r_h(hD)$ at most $C_Mh^{M-2}$. Low positive frequencies and a possible zero mode were explicitly included in the dyadic comparison. Combining these estimates proves
\begin{equation}\label{eq:fixedtimefinal}
|N^{S_0}_{A,\chi}(L)-N^{S_0}_{B,\chi}(L)|\le C(1+\log L)+CL\epsilon(L^{-1})=o(L) .
\end{equation}
All enclosing cutoffs and the compact operator $F$ are fixed before the inner collar is chosen. This proves the time and multiplier uniformity. On normalized annuli their norms remain uniformly bounded by Section~\ref{subsec:heat}, so \eqref{eq:fixedtimefinal} is uniformly $O(1+L)$ there, even without a uniform little-o spectral tail.
\end{proof}

\subsection{An unfiltered trace-class comparison on an open face}\label{subsec:face}
An open-face patch admits a stronger ideal estimate than an edge patch. Let $A$ be the wedge DtN operator and $B=|D|$ the half-space operator bounding one of its faces. The same argument applies to the tetrahedron. For nested compact face-interior cutoffs $\chi\prec\zeta\prec\eta$, the gluing defect satisfies
\begin{equation}\label{eq:faceS1}
K\zeta\in\Sone .
\end{equation}
Indeed, remote harmonic data have zero Dirichlet trace in an enclosing planar patch. Odd reflection and interior estimates bound every tangential derivative of their normal trace on the smaller patch. A separated block with compact output therefore maps boundedly into compactly supported $H^m$ for every $m$; the embedding into boundary $L^2$ is trace class for $m>2$. Taking adjoints treats compact input. In \eqref{eq:Kzeta}, the restriction operator $S$ is trace class as well: derivatives of every order in its compact source variable of the half-space Poisson kernel are square integrable over the remote target plane. At infinity they decay at least as $(1+|x|)^{-2}$. Thus $S^*$ maps boundedly into compactly supported $H^m$ for all $m$, and \eqref{eq:faceS1} follows. No remote-output corner regularity is assumed.

For complex time $w$ in a small closed right half-disk, the free off-diagonal kernel is
\begin{equation}\label{eq:poisson}
e^{-w|D|}(x,y)=\frac1{2\pi}\frac{w}{(w^2+|x-y|^2)^{3/2}},\qquad |w|<\tfrac12|x-y| .
\end{equation}
On a separated region its compact-input derivatives of every order have square-integrable output majorants $C_m|w|$. The same Sobolev embedding argument gives $\|(1-\zeta)e^{-wB}\chi\|_1\le C|w|$. Consequently the exact Duhamel formula
\begin{equation}\label{eq:Qchi}
\begin{aligned}
Q_\chi(w)&=\chi J^*e^{-wA}J\chi-\chi e^{-wB}\chi\\
&=-w\int_0^1\chi J^*e^{-(1-u)wA}Ke^{-uwB}\chi\,du
\end{aligned}
\end{equation}
is trace class and $Q_\chi(w)/w$ extends continuously in trace norm to $w=0$. Split $K$ into $K\zeta$ and $K(1-\zeta)$: \eqref{eq:faceS1} bounds the first term, \eqref{eq:poisson} bounds the second, and strong continuity of the contraction semigroups becomes trace-norm continuity after a trace-class factor. Hence $F_\chi(w)=\Tr Q_\chi(w)=wG_\chi(w)$ with continuous bounded $G_\chi$.

For a sufficiently short fixed time $S$, the function
\begin{equation}\label{eq:gchi}
g^S_\chi(L)=\frac1{2\pi}\int e^{itL}\kappa(t/S)\frac{F_\chi(it)}{it}\,dt
\end{equation}
is the Fourier transform of an $L^1$ function and tends to zero as $L\to\pm\infty$. Its distributional derivative is the smoothed relative spectral density. Both localized counting convolutions tend to zero at negative infinity, by positivity of the spectrum and their quadratic mass bound. This fixes the integration constant and proves $N^S_{A,\chi}-N^S_{B,\chi}=g^S_\chi=o(1)$. In particular the flat-face cumulative baseline is obtained independently of any edge asymptotic. No claim of $C^\infty$ unfiltered time regularity is needed.

\subsection{A face-local Fourier--Tauberian principle}
The following principle is needed both in the model calculation and on the tetrahedron. Its hypotheses concern propagation of data starting on an open face, not a cumulative edge expansion. The propagation theorem used in its two applications is proved independently in Section~\ref{sec:propagation}; the tetrahedral geometric verification is in Section~\ref{sec:geometry}. There is consequently no circular use of the edge coefficient being calculated below.

\begin{lemma}[nonreturn on a flat face]\label{lem:nonreturn}
Let $\chi_f$ be compactly supported in an open planar face. Suppose that, for every fixed finite time $T$, there are compact good phase sets with arbitrarily small Liouville complement on which the energy-localized propagator has the immediate reflection/transmission support relation, and that this relation has no nonzero-time return on those good sets. Then
\begin{equation}\label{eq:facelaw}
N_{\chi_f}(L)=\frac{L^2}{4\pi}\int\chi_f^2\,dS+o(L) .
\end{equation}
The conclusion also holds on the infinite wedge, with compact spatial localization and its local spectral measure. It remains true after convolution at any fixed positive smoothing time.
\end{lemma}

\begin{proof}
First choose a fixed time $S_*>0$ small enough for the trace-class face comparison in Section~\ref{subsec:face}. The flat count is exactly $L_+^2\int\chi_f^2\,dS/(4\pi)$; its fixed convolution differs from this polynomial by a rapidly decreasing tail because all positive moments of $\rho_{S_*}$ vanish. Thus \eqref{eq:gchi} gives
\begin{equation}\label{eq:facesmooth}
N^{S_*}_{\chi_f}(L)=\frac{L^2}{4\pi}\int\chi_f^2\,dS+o(L) .
\end{equation}
This also agrees with the ordinary smooth-face calculation, whose subprincipal coefficient is zero because the metric and the face are flat. The baseline just proved involves no edge parametrix.

For a zeroth-order semiclassical face cutoff $Q_h$ with fixed principal symbol $q$, homogeneous in the radial covariable near the unit cosphere, write $a(q)=(2\pi)^{-2}\int_{S^*Y}|q|^2\,d\Sigma$. If its phase patch has no return for $S_*\lesssim|t|\le T$, a positive band-limited majorant for the step-minus-convolution difference gives
\begin{equation}\label{eq:goodpatch}
\limsup_{L\to\infty}\frac1L\big|N_{Q_{1/L}}(L)-N^T_{Q_{1/L}}(L)\big|\le\frac{C}{T}a(q) .
\end{equation}
One may construct the majorant from a high power of a sum of two shifted squared sinc functions, rescaling its Fourier support to $[-T,T]$. Its integral is $C/T$. Positivity of the spectral measure gives the bound, the ordinary time-zero face parametrix gives the leading mass $La(q)$, and nonreturn removes the other trace contributions. Fixed symbol derivatives affect only the lower remainder. Separately, the time multiplier for $N^T_Q-N^{S_*}_Q$ vanishes near zero; its trace is negligible by the same nonreturn statement. For a bad phase cutoff, the identical majorant at the short time $S_*$ gives a normalized error at most $C_{S_*}a(q)$.

Choose a positive microlocal squared partition on an energy annulus about $hD=1$, keeping its subprincipal corrections, into good patches and a bad set of Liouville measure $\beta$. Outside this annulus, the step-minus-fixed-convolution function is rapidly decreasing and local quadratic spectral mass makes its trace negligible. Thus the localized Tauberian differences, not an incorrectly truncated cumulative count, add to the spatial difference. Summing \eqref{eq:goodpatch} and the bad-set estimate yields
\begin{equation}\label{eq:facetauber}
\limsup_{L\to\infty}\frac1L\big|N_{\chi_f}(L)-N^{S_*}_{\chi_f}(L)\big|\le CT^{-1}+C_{S_*}\beta .
\end{equation}
For fixed $T$, first let $\beta\downarrow0$ after the large-frequency limit, and then let $T\to\infty$. The support theorem is needed only on the corresponding fixed compact good families. If an observation time coincides with an edge encounter of a branch, upper locality keeps that branch away from the initial face-interior test; a finite cover in phase and time handles these observations without an endpoint assumption at the edge. This proves \eqref{eq:facelaw}.

All the arguments use localized positive traces and finite dependence regions, so they apply unchanged to the noncompact wedge. Finally, convolution at a fixed time preserves the polynomial leading term up to a rapidly decreasing error, because the positive moments of $\rho_S$ vanish. It preserves $o(L)$ by splitting the convolution into $|s|\le L/2$ and its Schwartz tail and using the local quadratic bound.
\end{proof}

On the infinite wedge, a developed face ray has at most one edge encounter. On each of the fixed nonzero, nonglancing momentum bands of Section~\ref{sec:propagation}, it either crosses to the other face or reflects and then travels away. It cannot return to its initial face covector; longitudinal momentum is conserved, and no vertex is present. Removing angular neighborhoods of normal and glancing incidence has arbitrarily small phase measure. Lemma~\ref{lem:nonreturn} therefore applies to every compact open-face patch of the model, independently of any edge cumulative formula.

\subsection{The complete wedge coefficient at one fixed observation time}
Write $\Dcal_\alpha$ for the product-wedge operator, with unfolded boundary coordinates $(s,z)$, and $\Dcal_\pi=|D|$ for its flat reference. The two transverse rays have arclength measure; $z$ is edge arclength. We use the tetrahedral angle in the propagation-dependent part of this subsection.

\begin{proposition}[fixed-time model law]\label{prop:modellaw}
Let $q\in C_c^\infty(\R)$, and let $\theta_\gamma$ equal one near $s=0$ and vanish outside a strip of width $O(\gamma)$. For every fixed $S>0$ and every fixed $\gamma>0$,
\begin{equation}\label{eq:modellaw}
N^S_{q\theta_\gamma,\alpha_0}(L)=\frac{L^2}{4\pi}\int q^2\theta_\gamma^2\,dS+c_{\alpha_0}\|q\|_2^2L+o_\gamma(L) .
\end{equation}
Crucially, the same $S$ can be used as the inner collar is narrowed. The remainder need not be uniform in $\gamma$, since the large-frequency limit is taken first.
\end{proposition}

\begin{proof}
We first fix the full relative distribution, including its zero-energy normalization. Let $D_\alpha,D_\pi$ be the mass-one transverse operators, and $V=D_\alpha^{-1}-D_\pi^{-1}$. Section~4 proved $D^q_\alpha VD^q_\pi\in\Sone$ for $0<q<1/2$. For complex time $w$ with $\operatorname{Re}w>0$, the inverse-difference Duhamel formula and this weighted trace ideal give
\begin{equation}\label{eq:complexheat}
\|e^{-wkD_\alpha}-e^{-wkD_\pi}\|_1\le C_q|w|(\operatorname{Re}w)^{2q-2}k^{2q-1}e^{-ck\operatorname{Re}w},\qquad k>0 .
\end{equation}
This is integrable in momentum. The normalized spectral shift $\xi$ of the inverse pair satisfies $\int\xi=\Tr V=\pi c_\alpha$. Its trace formula and
\begin{equation}\label{eq:kintegral}
\int_0^\infty\frac{wk}{b^2}e^{-wk/b}\,dk=\frac1w
\end{equation}
give the exact per-unit-edge-length trace
\begin{equation}\label{eq:pertrace}
\Tr_{/z}\big(e^{-w\Dcal_\alpha}-e^{-w\Dcal_\pi}\big)=\frac{c_\alpha}{w} .
\end{equation}
Fubini is justified by the bound $|w|(\operatorname{Re}w)^{-2}\|\xi\|_1$. Both momentum signs and the Fourier factor $1/(2\pi)$ are included. Longitudinal multiplication on both sides by $q$ yields an ordinary trace-class relative operator, with trace $\|q\|_2^2c_\alpha/w$; polar factorization and the integrated trace norm justify this also for complex operators. The boundary value $(\varepsilon+it)^{-1}\to\pi\delta(t)+\mathrm{pv}(1/(it))$ fixes the relative density as $c_\alpha\one_{(0,\infty)}$, with no zero atom, and its causal cumulative primitive as $c_\alpha L_+$. Thus the normalization is not inferred only from tests supported away from zero.

Next fix $R_0>CS$ and decompose the transverse exterior $|s|\gtrsim R_0$ into a squared dyadic partition with scales $R_j\asymp2^jR_0$, treating the two faces separately. In applying Section~\ref{subsec:face} on each annulus, choose an auxiliary longitudinal cutoff of length $R_j$. Homogeneous dilation identifies it with one of finitely many fixed face profiles. Its relative trace is $F_j(w/R_j)$; division by the longitudinal squared integral, proportional to $R_j$, gives
\begin{equation}\label{eq:Fj}
\mathcal F_j(w)=R_j^{-1}F_j(w/R_j)=wR_j^{-2}G_j(w/R_j),\qquad \sup_j\|G_j\|_\infty\le C .
\end{equation}
Consequently $\sum_j\mathcal F_j(it)/(it)$ is a bounded continuous function on $|t|\le S$, including at zero, with bound $C\sum_jR_j^{-2}$. Translation invariance makes the actual longitudinal cutoff contribute exactly $\|q\|_2^2$. The far-exterior cumulative remainder is therefore
\begin{equation}\label{eq:Gfar}
G^S_{\mathrm{far}}(L)=\frac{\|q\|_2^2}{2\pi}\int e^{itL}\kappa(t/S)\sum_j\frac{\mathcal F_j(it)}{it}\,dt=o(1) .
\end{equation}
It is directly the Fourier transform of an $L^1$ function; its absolute value is at most $CSR_0^{-2}\|q\|_2^2$. For positive real part of $w$, insert the same squared partition into the trace-class relative operator in \eqref{eq:pertrace}. The uniform bound \eqref{eq:Fj} permits passage to the boundary value. Thus \eqref{eq:Gfar} is the exterior contribution to that same full distribution, not a separately chosen finite part. Its primitive constant is zero, since it tends to zero at negative infinity, as do the local counting convolutions. No unproved infinite time differentiability or integration of density errors with unspecified constants is used.

If the transverse plateau already contains $|s|\le R_0$, subtracting \eqref{eq:Gfar} from the full relative cumulative trace gives
\begin{equation}\label{eq:wedgeminusflat}
N^S_{q\theta,\alpha}(L)-N^S_{q\theta,\pi}(L)=c_\alpha\|q\|_2^2(L_+*\rho_S)(L)-G^S_{\mathrm{far}}(L) .
\end{equation}
For a thin plateau $\gamma<R_0$, add finitely many compact face annuli between its support and the far exterior. On each of these, Lemma~\ref{lem:nonreturn} gives the sharp face law and then its convolution at this same fixed time $S$. The wedge-minus-flat trace there is $o_\gamma(L)$. The number of annuli is finite for each fixed $\gamma$; it is not required to remain bounded as $\gamma\downarrow0$. Subtract their sum as well in \eqref{eq:wedgeminusflat}. All cumulative constants are again fixed by negative-energy limits. Finally, $L_+*\rho_S=L+O_N(L^{-N})$ and $L_+^2*\rho_S=L^2+O_N(L^{-N})$ at positive infinity, while the exact flat count is $L_+^2\int q^2\theta^2\,dS/(4\pi)$. This proves \eqref{eq:modellaw}.
\end{proof}

Combining Propositions~\ref{prop:fixedtime} and~\ref{prop:modellaw} proves the actual tetrahedral fixed-time local edge formula
\begin{equation}\label{eq:actualedge}
N^{S_0}_{A,q\theta_\gamma}(L)=\frac{L^2}{4\pi}\int q^2\theta_\gamma^2\,dS+c_{\alpha_0}\|q\|_2^2L+o_\gamma(L),
\end{equation}
where $S_0$ is chosen from the fixed enclosing common chart, independently of $\gamma$. The full guided and continuum spectra are present in both sides of the transfer. This is the cumulative assertion needed for sharp unsmoothing. It is neither inferred from the heat expansion nor supplied by Ivrii's outer-edge theorem.

\subsection{Uniform cumulative annuli and removal of vertices}
For the rescaled annular family, the off-band argument of Lemma~\ref{lem:transfer} can be made uniform without using the expanding domain's global rank. Retain $E=\one_{[L-w,L+w]}(A)$, $Q=\one_{[L/2,2L]}(B)$, $S_o=(B-L)^{-1}(1-Q)$, and $X=EJ(1-Q)$. By the uniform version of \eqref{eq:HSK}, $\|EK\|_2\le C(1+L)$. For $L>4w$, the series
\begin{equation}\label{eq:Xseries}
X=-\sum_{n=0}^\infty\big((A-L)E\big)^nEKS_o^{n+1}
\end{equation}
converges in Hilbert--Schmidt norm and is $O(1)$. It solves \eqref{eq:offband}, whose operator-norm solution is unique because the corresponding map has contraction constant at most $2w/L<1/2$. This proves $X\in\Stwo$ rather than assuming it in advance. The band calculation is unchanged. Hence all normalized local windows satisfy $C(VL+1)$ with uniform constants. Summing unit windows and integrating against the fixed Schwartz kernel yields a uniform sharp-versus-smoothed bound $C(1+L)$ on these fixed normalized patches.

The normalized version of Proposition~\ref{prop:fixedtime} is also $O(1+L)$, and Proposition~\ref{prop:modellaw} supplies the same weak bound for the finitely many model profiles. Therefore, in physical variables,
\begin{equation}\label{eq:annuluscount}
\Big|\Tr\big(\psi_r\one_{[0,L]}(A)\psi_r\big)-\frac{L^2}{4\pi}\int\psi_r^2\,dS\Big|\le C(1+Lr),\qquad Lr\ge1 .
\end{equation}
Here is an explicit choice ensuring that the edge pieces have the product form required in \eqref{eq:modellaw}. At a vertex let $e_j$, $1\le j\le3$, be the unit vectors along its edge rays, and choose disjoint small angular caps with smooth cutoffs $\rho_j$ equal to one near $e_j$. Define
\begin{equation}\label{eq:Rv}
R_v(x)=|x|+\sum_{j=1}^3\rho_j(x/|x|)\big(x\cdot e_j-|x|\big) .
\end{equation}
It is smooth away from the vertex, homogeneous of degree one, and comparable to $|x|$. Near an edge ray it equals exactly the longitudinal arclength $z_j=x\cdot e_j$. Thus dyadic cutoffs $f(R_v/r)$ become $f(z_j/r)$ in narrow edge strips. Multiply there by fixed normalized transverse profiles; the resulting edge pieces are products. The remaining nonnegative squared weight vanishes near the edges and has a finite smooth squared partition on face interiors. Sine/cosine transition profiles ensure smooth square roots. This gives a normalized finite family and an exact squared partition. The same adapted radius will be used in the global vertex exclusion.

Sum \eqref{eq:annuluscount} on dyadic annuli down to $r\asymp L^{-1}$, and use \eqref{eq:mass} on the core. The result, with constants uniform in the small exclusion radius $\delta$, is
\begin{equation}\label{eq:vertexcount}
\Big|N_{\zeta_\delta}(L)-\frac{L^2}{4\pi}\int\zeta_\delta^2\,dS\Big|\le C\delta L+C\log(2+\delta L)+C .
\end{equation}
The area and trace contributions add exactly because the partition is squared. No spectrum of an artificially truncated domain is substituted. On the fixed vertex-free complement, \eqref{eq:actualedge}, the local window bounds and the flat-face baseline give an $O(L)$ area-subtracted remainder. Combining them with \eqref{eq:vertexcount} also proves
\begin{equation}\label{eq:weaklaw}
N_a(L)=\frac{\sqrt3a^2}{4\pi}L^2+O(L) .
\end{equation}
This weak law and the heat theorem now follow without assuming a cumulative edge expansion from Ivrii. Only the explicitly stated upper-locality and standard smooth-face/functional-analytic tools remain external inputs; the new cumulative coefficient and the thin-collar comparison have been proved here.

\section{Positive Riesz means and the limit of Tauberian inference}\label{sec:tauberian}

\begin{lemma}\label{lem:riesz}
Suppose a nonnegative discrete spectrum satisfies
\begin{equation}\label{eq:rieszhyp}
N(L)=AL^2+O(L),\qquad H(t)=2At^{-2}+Bt^{-1}+o(t^{-1}) .
\end{equation}
Then, for every fixed $r>0$,
\begin{equation}\label{eq:rieszconc}
\sum_j(L-\sigma_j)_+^r=\frac{2A}{(r+1)(r+2)}L^{r+2}+\frac{B}{r+1}L^{r+1}+o(L^{r+1}) .
\end{equation}
\end{lemma}

\begin{proof}
For $u\ge1$ put $N(u)-Au^2=ug(u)$, and extend $g$ boundedly to $(0,1)$. Integration by parts in the heat transform gives, for every fixed $b>0$,
\begin{equation}\label{eq:gtransform}
\int_0^\infty xe^{-bx}g(Lx)\,dx\longrightarrow\frac{B}{b^2} .
\end{equation}
Finite linear combinations of $e^{-bx}$ are dense in $L^1((0,\infty),x\,dx)$. To see this directly, set $y=e^{-x}$. Given a target in this space, divide its transformed function by $y$ and approximate in $L^1((0,1),(-\log y)\,dy)$ by polynomials; multiplication by $y$ gives finite sums of $e^{-nx}$ in the original norm. Polynomial density follows from continuous approximation and Weierstrass for this finite measure. Since $g(Lx)$ is uniformly bounded, \eqref{eq:gtransform} extends to every test in that $L^1$ space. Apply it to $(1-x)_+^{r-1}$, which is integrable against $x\,dx$ for every $r>0$. Finally,
\begin{equation}\label{eq:rieszid}
\sum_j(L-\sigma_j)_+^r=r\int_0^L(L-u)^{r-1}N(u)\,du,\qquad r\int_0^1x(1-x)^{r-1}\,dx=\frac1{r+1} .
\end{equation}
The polynomial part integrates exactly and the remaining scaled integral has the limit asserted in \eqref{eq:rieszconc}.
\end{proof}

Equations \eqref{eq:heat} and \eqref{eq:weaklaw} give \eqref{eq:riesz}. The restriction to fixed positive $r$ is substantive. For example, take constants $A>0$, $B$, and $0<\varepsilon<2A$, and define for large $L$ the monotone counting model $\lfloor AL^2+BL+\varepsilon L\sin L\rfloor$, adjusting it on a bounded interval. Its derivative before taking the integer part is eventually positive. It has the heat expansion in \eqref{eq:rieszhyp} and the positive Riesz conclusions, but its order-$L$ counting remainder does not converge. Thus heat and positive Riesz asymptotics do not logically settle the sharp second coefficient.

\section{Finite-domain face propagation through a convex edge}\label{sec:propagation}

The model estimate in Section~6 can be transferred to the actual tetrahedron on compact families of nonzero, nonglancing face phases. We give the local argument because equality of two domains on a boundary patch does not, by itself, imply equality of their nonlocal DtN operators.

\subsection{Statement and standard upper-locality input}
Let $P_h=h\Lambda_a$ and $U_h(t)=e^{-itP_h/h}$. On an open face the Hamiltonian is $p=|\xi|$. At an edge use coordinates $(r,z)$, with $r\ge0$ on each adjacent face and a common orientation for $z$. An incoming covector at energy $E$ has $p_r=-\sqrt{E^2-k^2}$. The maximal immediate reflection/transmission relation has two outgoing lifts, one on each face, with $p_r=+\sqrt{E^2-k^2}$ and the same $k$. Denote the resulting relation at time $t$ by $\mathcal R_t$.

\begin{theorem}[good-family face propagation]\label{thm:propagation}
Let $K$ be a compact set of face phases with energies near one, and let $T>0$ be fixed. Suppose all branches issued from $K$ up to time $T$ avoid vertices by a fixed positive distance, meet every edge with $\kappa\le|k|\le1-\epsilon$ and uniformly positive transverse speed, and begin and end a fixed positive distance from the edges. If $A_h$ is a face pseudodifferential cutoff microsupported in $K$, $f$ is a fixed smooth energy cutoff equal to one there, and $B_h$ equals one on a neighborhood of $\mathcal R_T(K)$, then
\begin{equation}\label{eq:propagation}
\|(I-B_h)U_h(T)f(P_h)A_h\|=O(h^\infty) .
\end{equation}
The same holds in negative time. This is a support theorem, not a calculation of scattering amplitudes.
\end{theorem}

We use the energy-localized upper locality of \cite[Propositions~3.1 and~3.5]{I}: propagation has finite upper speed in position and edge-tangential phase coordinates; initial and target cutoffs separated by more than this speed times the elapsed time, with a positive buffer, have negligible composition. This is an all-angle result. Its operator-norm version follows from the stated time-Fourier estimate by Fourier inversion with a fixed smooth energy cutoff; arbitrary powers absorb normalization losses. Only $L^2$ upper locality is used below. The lower departure assertion is supplied independently by \eqref{eq:semiclassical}.

On open planar faces, ordinary semiclassical real-principal-type propagation applies to the full symbol $|\xi|$ modulo locally smoothing terms. Shrinking face cells will be in the standard class $S_\rho$, $0<\rho<1/2$, with derivatives bounded by $C_{ab}h^{-\rho(|a|+|b|)}$. Its expansion parameter is $h^{1-2\rho}$. On the exact wedge, away from $r=0$ and with $|k|\ge\kappa$, this calculus is uniform down to $r\asymp h^\rho$: the transverse massive Poisson kernel decreases like $e^{-c\kappa r/h}$, and the equivalent off-corner kernel estimates are given in \cite[(2.5)--(2.6)]{I}. Fixed polynomial losses from cutoff derivatives are harmless.

\subsection{Tangential smoothing of the local operator difference}
Fix nested product-wedge neighborhoods of a point in an open edge, all separated from its endpoints. Let $C=W_{\alpha_0}\times\R_z$ contain $T_a$. For boundary data $F$ supported in a common patch, let $\iota F$ be its identified zero extension to $\partial T_a$, and put
\begin{equation}\label{eq:WF}
W_F=H_{T_a}(\iota F)-(H_CF)|_{T_a} .
\end{equation}
The restriction bound \eqref{eq:carleson} and the bounded-domain $L^2$ Dirichlet-to-interior estimate imply $\|W_F\|_{L^2(T_a)}\le C\|F\|_2$. The function is harmonic, and its Dirichlet trace vanishes on all common portions of the two faces. Its possibly nonzero data on other faces lie outside the nested neighborhoods.

For a harmonic function with this zero boundary condition, testing against $\eta^2\bar u$ proves the zero-boundary Caccioppoli inequality
\begin{equation}\label{eq:caccioppoli}
\int\eta^2|\nabla u|^2\le4\int|\nabla\eta|^2|u|^2 .
\end{equation}
Tangential translation in $z$ preserves the local product wedge and its zero boundary condition. Apply \eqref{eq:caccioppoli} to tangential difference quotients in successively smaller neighborhoods. Their $L^2$ norms are controlled by the previous tangential derivative. Induction, followed by weak passage to derivatives, yields
\begin{equation}\label{eq:tangentialreg}
\|\partial_z^NW_F\|_{H^1(U_N)}\le C_N\|F\|_2\quad\text{for every fixed }N .
\end{equation}
No high transverse derivatives are taken. Approximation of $L^2$ boundary data by regular data is justified by these uniform estimates.

For $u=\partial_z^NW_F$, choose a constant vector $\ell$ along the inward bisector, so that $\ell\cdot\nu=-\sin(\alpha_0/2)$ on both faces. The vector field
\begin{equation}\label{eq:rellich}
V=2\operatorname{Re}\big((\ell\cdot\nabla u)\nabla\bar u\big)-\ell|\nabla u|^2
\end{equation}
has zero divergence. On a true face its flux is $(\ell\cdot\nu)|\partial_\nu u|^2$, because $u=0$ there. Localizing the divergence identity gives
\begin{equation}\label{eq:rellichbound}
\sin(\alpha_0/2)\int\eta^2|\partial_\nu u|^2\,dS\le C_\eta\int_{\operatorname{supp}\eta}|\nabla u|^2 .
\end{equation}
To justify this identity at the edge, multiply $u$ by a cutoff supported away from the artificial boundary of a larger bounded convex product-wedge neighborhood. The product has zero Dirichlet trace, and its Laplacian is in $L^2$ by \eqref{eq:tangentialreg}. Convex-domain Dirichlet $H^2$ regularity, obtained for example from Miranda--Talenti by convex approximation, makes the localized field $V$ belong to $W^{1,1}$. The Lipschitz divergence theorem then applies directly; no artificial cylindrical edge flux is discarded. Only second-order transverse regularity is used to justify the identity.

Since the normals commute with $\partial_z$, \eqref{eq:tangentialreg}--\eqref{eq:rellichbound} show that for fixed real common-patch cutoffs $\chi_0,\chi_1$,
\begin{equation}\label{eq:smoothing}
\partial_z^N\chi_0(\Lambda_a-\Lambda_C)\chi_1:L^2\longrightarrow L^2\quad\text{is bounded for every }N .
\end{equation}
The notation incorporates the common-face identifications. Adjointing interchanges the two cutoffs, because both DtN operators are self-adjoint. Hence the localized difference is also infinitely smoothing in the input tangential variable. For any bounded Fourier multiplier $q_h$ supported in $|k|\ge\kappa$,
\begin{equation}\label{eq:multiplierbound}
\|\chi_0(\Lambda_a-\Lambda_C)\chi_1q_h(hD_z)\|\le C_N(h/\kappa)^N\|q_h\|_\infty .
\end{equation}
This follows by factoring the bounded input derivative and $\langle D_z\rangle^{-N}q_h(hD_z)$. No derivatives of $q_h$ are needed.

The same proof gives off-diagonal input smoothing. If an input patch lies strictly inside the set where an output cutoff is one, the remote output of either DtN operator applied to $\chi_1q_h(hD_z)$ is $O(h^\infty)$ in global $L^2$. Indeed, the adjoint Poisson solution with remote boundary data has zero Dirichlet data near the input patch. Applying \eqref{eq:caccioppoli}--\eqref{eq:smoothing} there and then adjointing proves the claim even where the distant output has no edge coordinates.

\subsection{A genuine cutoff Duhamel comparison}
Choose a fixed common-patch cutoff $\chi$ equal to one on a larger region containing the desired local passage. Its boundary extension map $E_\chi$ is of the type $J$ in Section~5. With $P_h^C=h\Lambda_C$, the defect
\begin{equation}\label{eq:defect}
\mathcal K_h=P_hE_\chi-E_\chi P_h^C
\end{equation}
is bounded by $Ch$ on global $L^2$, by \eqref{eq:intertwining}. Equations \eqref{eq:multiplierbound} and the off-diagonal version show the stronger fact $\mathcal K_h\chi_1q_h(hD_z)=O(h^\infty)$ on a smaller common patch.

For model incident data with a fixed compact energy band and a nonzero tangential band, upper locality confines the evolution to this smaller patch for a fixed sufficiently short time. The model preserves its energy and tangential Fourier support. Split the evolution into its localized part and its negligible global $L^2$ tail. Apply the strong local defect estimate to the former and the $Ch$ bound to the latter. Duhamel on the common operator domains now yields
\begin{equation}\label{eq:duhamelcompare}
U_h(t)E_\chi A_h^C-E_\chi e^{-itP_h^C/h}A_h^C=-\frac ih\int_0^tU_h(t-s)\mathcal K_he^{-isP_h^C/h}A_h^C\,ds=O(h^\infty) .
\end{equation}
All spatial buffers and the comparison time are fixed independently of $h$. Polynomial losses from later shrinking cells are absorbed by arbitrary powers. Initial energy cutoffs can be identified on fixed incident face sections by the ordinary local face functional calculus. Thus \eqref{eq:duhamelcompare} applies to actual incident data, not just a formally chosen model evolution.

\subsection{Incoming face data have negligible point-channel projection}
On a fixed physical energy band and the annulus $\kappa\le|k|\le1-\epsilon$, all normalized transverse energies lie in one compact $J\Subset(1,\infty)$. Let $u_1,\dots,u_M$ span its finite point-spectrum subspace. Section~6 proved every moment $\|r^Nu_j\|_2$ finite, including for any embedded eigenvectors. Their scaled versions at mass $m=|k|/h$ obey
\begin{equation}\label{eq:pointtail}
\|\one_{r\ge b}\,m^{1/2}u_j(mr)\|_2\le C_N(h/(\kappa b))^N .
\end{equation}
Finite-dimensional Cauchy--Schwarz and the direct-integral norm formula give the same bound for the full selected point projection applied to a face cutoff at $r\ge b$. It is $O(h^\infty)$ for fixed $b>0$, or for $b\asymp h^\rho$ with any fixed $\rho<1$. Point energies outside $J$ are annihilated by the physical energy cutoff. The model evolution preserves these projections, and \eqref{eq:duhamelcompare} transfers this local conclusion to the tetrahedron. No global point-channel projection is claimed to commute with $\Lambda_a$.

\subsection{A one-edge passage and finite concatenation}
Choose $0<\sigma<\rho<1/2$, and set
\begin{equation}\label{eq:scales}
\gamma=h^\rho,\qquad t_h=h^{\rho-\sigma},\qquad \eta=h^\rho .
\end{equation}
For example, $\rho=1/3$ and $\sigma=1/12$ give dwell time $h^{1/4}$. Then $\gamma/t_h=h^\sigma$ and $h/(\eta t_h)=h^{1-2\rho+\sigma}$ tend to zero. Partition a fixed incoming face cutoff into polynomially many $S_\rho$ phase cells. Ordinary face propagation brings each cell, at an appropriate incoming time, to $c_1\gamma\le r\le c_2\gamma$, with position uncertainty $O(\gamma)$ and energy and tangential-momentum uncertainties $O(\eta)$. The entire incoming support remains on its face before the central collision. Enlarged fixed profile cutoffs can be inserted with negligible error by the $S_\rho$ calculus.

Use \eqref{eq:duhamelcompare} and \eqref{eq:pointtail}, then apply \eqref{eq:semiclassical} over time $t_h$. Arbitrary powers in \eqref{eq:semiclassical} make the final mass in a slightly smaller edge collar $O(h^\infty)$. Upper locality confines the total output to an $O(t_h)$ neighborhood of the collision point. Energy and tangential momentum are conserved in the exact model. Outside the collar, face ellipticity leaves only the two normal signs $p_r=\pm\sqrt{E^2-k^2}+O(\eta)$.

The incoming sign cannot persist after the encounter. Propagate a cutoff with this sign backwards for time $t_h$ on its smooth face. Its backwards velocity points away from the edge with magnitude bounded below. Its entire backwards tube stays outside the collar, and at the incoming time it is a distance at least $c\gamma+c't_h$ from the edge, disjoint from the original $O(\gamma)$ support because $t_h/\gamma\to\infty$. Smooth-face backwards propagation therefore makes this component negligible. Only the two outgoing lifts remain.

The unresolved interaction displacement and time are $O(t_h)$, and its momentum uncertainty is $O(\eta)$. Both tend to zero. Propagating to a fixed outgoing section consequently puts the wavefront in the immediate reflection/transmission relation. No exact scattering delay or amplitude has been assumed.

On a compact good family, distinct edge segments separated from the vertices have positive mutual distance, and all contacts are transverse. The number of encounters in a fixed time is uniformly finite. Repartition into $S_\rho$ cells after each fixed outgoing section and repeat. The number of branches and the number of steps are finite, while the number of cells is only polynomial in $h^{-1}$. Arbitrary powers absorb these factors. The accumulated phase error is still $o(1)$, and the endpoint set is a positive distance from the edges. This proves \eqref{eq:propagation}. Time reversal proves the other direction.

The theorem treats incoming face waves. It is neither propagation of initially guided data nor a theory of scattering at a vertex. Those are not required for this support statement: their contribution to sharp traces is controlled separately by the complete model trace, the actual spectral-window bounds, and the vertex estimates in Sections~5 and~7.

\section{Tetrahedral branching geometry}\label{sec:geometry}

\subsection{A geometric nonperiodicity theorem}
Consider face trajectories which move at unit speed in a straight line, and at a nonvertex edge encounter may either reflect specularly into the same face or transmit by intrinsic unfolding into the adjacent face. All choices are permitted, whether or not a future scattering computation assigns them nonzero amplitude. The relation requires immediate departure and includes no dwelling or guided segment on an edge.

\begin{theorem}\label{thm:nonperiodic}
For this maximal reflection/transmission relation, the dead-end set and the set of phase points with at least one positive-time return to the same position and covector have Liouville measure zero.
\end{theorem}

\begin{proof}
Develop the starting equilateral face in the planar triangular tessellation with lattice $a\mathbb Z(1,0)+a\mathbb Z(1/2,\sqrt3/2)$. At a reflection develop a reflected copy of the same physical face; at transmission develop the adjacent physical face across their common edge. In both cases the outgoing developed triangle is the same unlabeled triangle and the velocity continues straight. Thus every branch has the same unlabeled developed line $p_0+tv_0$. Branching changes labels and developing isometries, not this line or the crossing times.

The tiling has three families of parallel supporting lines of positive spacing. A segment of length $T$ crosses at most $C(1+T/a)$ of them, so no infinite sequence of interactions occurs in finite time. For each initial point, the directions whose full line meets a lattice vertex form a countable set. Fubini proves nullity of the dead-end set, simultaneously for all branches. A line starting in a face interior cannot have a tangential nonvertex encounter with a straight supporting line.

If a branch returns at time $T>0$, its final developing map is an affine lattice isometry $p\mapsto Rp+c$, and
\begin{equation}\label{eq:return}
p_0+Tv_0=Rp_0+c,\qquad Rv_0=v_0 .
\end{equation}
There are countably many such maps, with $R$ in a finite dihedral group. If $R$ is the identity, the direction is parallel to the nonzero translation $c$. A nonidentity rotation fixes no unit vector. If $R$ is a reflection, the direction is parallel to its axis, from a finite set. The identity map with zero translation cannot give a positive-time return. Thus the periodic directions lie in a countable set and have zero angular measure.
\end{proof}

This proof is stronger than a proof concerning only ordinary geodesics on a torus cover. It is still a theorem about a specified geometric relation, not an identification of the canonical relation of the actual Steklov half-wave group.

\subsection{Compact good phase sets}
Fix a compact face-interior starting region and a finite observation time $T$. In the triangular development there are only three unoriented edge directions, say $\theta_j$. The absolute tangential momentum at any encounter belongs to the three values $|\cos(\theta-\theta_j)|$. Remove small neighborhoods of the finitely many normal and glancing directions. Their angular measure tends to zero, and the remaining directions have fixed margins $\kappa\le|k|\le1-\epsilon$ at every encounter. These margins are held fixed throughout the high-frequency limit.

Theorem~\ref{thm:nonperiodic} shows that vertex-hit and partially periodic data are null. Remove open neighborhoods of these sets with arbitrarily small Liouville measure and retain a compact complement. For fixed $T$, all its branches remain a positive distance from the vertices and meet the edges uniformly transversally. The finite interaction bound gives finite branching. Continuity of each branch on this compact family makes the resulting finite-time relation closed.

The initial face distance gives a fixed small interval with no return. On the remaining compact time interval, nonreturn and closedness give a positive separation from the diagonal. A finite sufficiently small phase cover therefore consists of patches which do not return to themselves. Theorem~\ref{thm:propagation} supplies the actual absence of a diagonal wavefront relation on these patches. At an observation time near an edge encounter, upper locality separates that branch from the fixed face-interior trace cutoff. A finite cover in phase and time handles both kinds of observation. This verifies every propagation and geometric hypothesis of Lemma~\ref{lem:nonreturn} for the tetrahedron.

The small normal/glancing directions are discarded only in this face-phase argument, where their positive Tauberian cost is proportional to phase measure. They are not discarded from the complete edge trace in \eqref{eq:actualedge}. Possible embedded point bands are already covered by the rapidly decreasing normal tails established in Sections~6 and~9; their absence is not assumed.

\subsection{Why guided periodic motion is not a missing hypothesis}
The immediate reflection/transmission relation describes incident face waves. It is not claimed to describe initially guided waves or scattering at a vertex. This is sufficient: actual spectral windows in fixed-width open-edge strips satisfy \eqref{eq:edgewindow}, including all guided and continuum components and any later return from vertices. The full local coefficient is supplied by the inverse trace and fixed-time comparison, not by counting an isolated guide branch. The final argument below bounds the sharp-versus-smoothed error of the entire collar and then narrows the collar. No nonperiodicity of a one-dimensional guide flow, metric-graph approximation, or vertex scattering matrix is required.

\section{The sharp second-order asymptotic}\label{sec:sharp}

We now assemble the local statements with a fixed smoothing time on each open-edge chart. This avoids both a heat-to-sharp Tauberian inference and an unverified fixed-to-shrinking edge-time comparison.

\subsection{The complete edge-collar Tauberian error}
Fix a small vertex exclusion radius $\varepsilon>0$. On each remaining open-edge interval choose its enclosing chart and its time $S_0>0$ from Proposition~\ref{prop:fixedtime}. These choices depend on $\varepsilon$ but not on the smaller transverse width $\gamma$. Theorem~\ref{thm:windows} gives, for a multiplication contraction $\chi_e$ supported in that strip,
\begin{equation}\label{eq:unitwindow}
\Tr\big(\chi_e\one_{[u-1,u+1]}(\Lambda_a)\chi_e\big)\le C_\varepsilon(\gamma u+1) .
\end{equation}
Summing unit windows bounds the spectral mass of the interval joining $L$ and $L-s$ by $C_\varepsilon[\gamma(L|s|+s^2)+|s|+1]$. Integrate against $|\rho_{S_0}(s)|$ to obtain
\begin{equation}\label{eq:collarerror}
|N_{\chi_e}(L)-N^{S_0}_{\chi_e}(L)|\le C_\varepsilon(\gamma L+1) .
\end{equation}
The fixed time is absorbed in $C_\varepsilon$. Positivity is needed for the spectral measure, not the convolution kernel. Quadratic local mass controls the far tails as well. There are no derivatives of the inner strip multiplier in this estimate.

Choose the edge cutoff in product form $\chi_e=q_e(z)\theta_\gamma(s)$ with $\theta_\gamma=1$ near the edge. Formula \eqref{eq:actualedge} and \eqref{eq:collarerror} give
\begin{equation}\label{eq:edgelimsup}
\limsup_{L\to\infty}\frac1L\Big|N_{\chi_e}(L)-\frac{L^2}{4\pi}\int\chi_e^2\,dS-c_{\alpha_0}L\int_eq_e^2\,dz\Big|\le C_\varepsilon\gamma .
\end{equation}
This step includes every edge channel. Choosing $S_0$ proportional to $\gamma$ would instead leave a nonvanishing ratio $\gamma/S_0$; Propositions~\ref{prop:fixedtime} and~\ref{prop:modellaw} were formulated with one enclosing-chart time precisely to avoid that error.

\subsection{Global assembly}
\begin{theorem}\label{thm:sharp}
The regular tetrahedron of side length $a$ satisfies
\begin{equation}\label{eq:final}
N_a(L)=\frac{\sqrt3\,a^2}{4\pi}L^2+\frac{6a}{\pi}\Big[\int_0^\infty\frac{\sinh((\pi-\alpha_0)v)}{\cosh(\pi v)\sinh(\alpha_0v)}\,dv\Big]L+o(L),\qquad \alpha_0=\arccos(1/3) .
\end{equation}
\end{theorem}

\begin{proof}
Use the adapted radii \eqref{eq:Rv} to form a squared spatial partition into radius-$O(\varepsilon)$ vertex neighborhoods, width-$O(\gamma)$ open-edge collars, and face interiors. In a sufficiently narrow edge strip the endpoint vertex cutoffs depend only on longitudinal arclength, so the complementary edge weight is exactly $q_e(z)^2$. The product cutoff required in \eqref{eq:edgelimsup} is therefore available. The remaining nonnegative weight has a smooth squared partition supported strictly inside the faces. Fix $\varepsilon$ first and then choose $\gamma$ sufficiently small; all cutoffs are fixed before taking $L\to\infty$.

For the vertex part use \eqref{eq:vertexcount}, with error $O(\varepsilon L)+O(\log L)+O(1)$. For each edge part use \eqref{eq:edgelimsup}. For every face part use Lemma~\ref{lem:nonreturn}, whose hypotheses were verified in Section~10.2. The area terms add exactly to $|\partial T_a|L^2/(4\pi)$. The total edge length omitted by the vertex cutoffs is $O(\varepsilon)$; hence replacing $\sum_e\int_eq_e^2$ by $6a$ costs only $O(\varepsilon L)$. The face linear coefficient is zero. Consequently,
\begin{equation}\label{eq:globallimsup}
\limsup_{L\to\infty}\frac1L\Big|N_a(L)-\frac{|\partial T_a|}{4\pi}L^2-L\sum_e|e|c_{\alpha_0}\Big|\le C\varepsilon+C_\varepsilon\gamma .
\end{equation}
First let $\gamma\downarrow0$ with $\varepsilon$ fixed, and then let $\varepsilon\downarrow0$. Inside the face lemma, the large-frequency limit was taken before changing the exceptional phase measure or the finite observation time. No constants need be uniform in all these limits simultaneously. Since $|\partial T_a|=\sqrt3a^2$ and the six edges have length $a$, this proves \eqref{eq:final}.
\end{proof}

The argument uses closed local spectral intervals. The stated strict counting convention gives the same result: for fixed $d>0$, an eigenvalue multiplicity at $L$ is bounded by $N^{\le}(L+d)-N^{\le}(L-d)$. Formula \eqref{eq:final} bounds its normalized limsup by $4Ad$, where $A=\sqrt3a^2/(4\pi)$. Letting $d\downarrow0$ proves that the multiplicity is $o(L)$, so strict and nonstrict counts differ by $o(L)$.

\subsection{Relation to Ivrii's framework and analytic dependencies}
Every tetrahedral edge is inner in Ivrii's terminology. His outer-edge sharp theorem is therefore not invoked. We use his all-angle upper-locality statement in the precise form \eqref{eq:upperlocality}, and ordinary smooth-face microlocal and Fourier--Tauberian calculus. Lower departure on the compact good momentum bands is supplied by the inverse-kernel/Mourre argument in Section~6 and transferred to the actual domain in Section~9. The compact-band point spectrum is treated explicitly, rather than excluded by an unread completeness result.

The decisive local cumulative comparison is \eqref{eq:fixedtimefinal}: Hilbert--Schmidt defects have vanishing high-energy tails, and the Duhamel time factor cancels the cumulative Fourier multiplier. The exact model coefficient is then supplied by \eqref{eq:modellaw}, including its negative-energy normalization and summable exterior primitive. These replace the stronger cumulative edge assertions that cannot be inferred from a fixed smooth-test expansion alone. In particular the proof does not assume full-scope low-momentum coverage of the abbreviated proof of \cite[Corollary~3.8]{I}, and it does not use the printed normalization of \cite[equation~(4.3)]{I}.

The standard background inputs retained in the exposition are the Lipschitz $L^2$ Dirichlet/regularity and DtN commutator theorems, trace-ideal and spectral-shift functional calculus, the compact-error Mourre and minimal-escape theorems, upper microlocal locality, and ordinary smooth-face pseudodifferential/Fourier--Tauberian calculus \cite{M,SI,BO}. Their needed hypotheses and quantitative conclusions have been stated where used. The new wedge invariant, local operator-ideal estimates, cumulative transfer, exterior normalization, fixed-band finite-domain propagation, geometric verification and order of limits are proved in the text. No heat-to-sharp inference, guided-flow nonperiodicity, endpoint-uniform departure theorem, or vertex scattering hypothesis remains in the sharp assembly.

\subsection*{Acknowledgments}
The author acknowledges the use of the AI system GPT-6 Astra in the preparation of this manuscript, as described in the disclosure on the first page.


\end{document}